\documentclass[12pt]{article}

\usepackage{amsfonts}
\usepackage{latexsym}
\usepackage[dvips]{graphicx} 
\usepackage[dvips]{color} 

\def\medskipamount{8pt} 
\def\arraystretch{1.3}
\begin{document}


\renewcommand{\thesection}%
{Part \Roman{section}\hspace{.2ex}:\hspace{-1ex}}
\renewcommand{\thesubsection}{\arabic{subb}}

\newcounter{subb}
\setcounter{subb}{0}

\newcommand{\bit}[1]{\pagebreak[3] \addtocounter{subb}{1}
  \subsection{\hspace{-2ex} #1}\setcounter{equation}{0}\penalty12000}
\renewcommand{\theequation}{\thesubsection.\arabic{equation}}


\newcommand{\re}[1]{\mbox{\bf (\ref{#1})}}


\catcode`\@=\active
\catcode`\@=11
\def\@eqnnum{\hbox to .01pt{}\rlap{\bf \hskip -\displaywidth(\theequation)}}
\catcode`\@=12



\catcode`\@=\active
\catcode`\@=11
\newcommand{\nc}{\newcommand}

   
\nc{\bs}[1]{\addpenalty{-1000}
\addvspace{\medskipamount}\penalty10000 
\refstepcounter{equation}
\noindent \begin{em}{\bf (\theequation) #1.} \nopagebreak}

\nc{\es}{\par \end{em} \addvspace{\medskipamount} } 

\nc{\br}[1]{\addvspace{\medskipamount} 
\refstepcounter{equation} 
\noindent {\bf (\theequation) #1.} \nopagebreak }

\nc{\er}{\par \addvspace{\medskipamount} }

\newcounter{index}
\nc{\bl}{\begin{list}{{\rm (\roman{index})}}{\usecounter{index}}}
   
\nc{\el}{\end{list}}

\nc{\pf}{\addvspace{\medskipamount} \par \noindent {\em Proof}}

\nc{\fp}{\phantom{.} \hfill \mbox{     $\Box$} \par
\addvspace{\medskipamount}} 

\nc{\beq}{\begin{equation}}
\nc{\eeq}{\end{equation}}

\nc{\beqas}{\begin{eqnarray*}}
\nc{\eeqas}{\end{eqnarray*}}


\nc{\C}{\mathbb C}
\nc{\Hyp}{\mathbf H}
\renewcommand{\P}{\mathbb P}
\nc{\Q}{\mathbb Q}
\nc{\R}{\mathbf R}
\nc{\Z}{\mathbb Z}
 

\nc{\op}[1]{\mathop{\mathchoice{\mbox{\rm #1}}{\mbox{\rm #1}}
{\mbox{\rm \scriptsize #1}}{\mbox{\rm \tiny #1}}}\nolimits}
\nc{\Ad}{\op{Ad}}
\nc{\ad}{\op{ad}}
\nc{\const}{\op{const.}}
\nc{\ch}{\op{ch}}
\nc{\pdeg}{\op{pdeg}}
\nc{\End}{\op{End}}
\nc{\END}{\op{$\mathbf{End}$}}
\nc{\Ext}{\op{Ext}}
\nc{\et}{\op{\'et}}
\nc{\Hom}{\op{Hom}}
\nc{\HOM}{\op{$\mathbf{Hom}$}}
\renewcommand{\Im}{\op{Im}}
\nc{\id}{\op{id}}
\nc{\Jac}{\op{Jac}}
\nc{\Map}{\op{Map}}
\nc{\mod}{\, / \,}
\nc{\Nm}{\op{Nm}}
\nc{\Pic}{\op{Pic}}
\nc{\Pol}{\op{Pol}}
\nc{\Quot}{\op{Quot}}
\nc{\slope}{\op{slope}}
\nc{\Triv}{\op{Triv}}
\nc{\rk}{\op{rk}}
\nc{\td}{\op{td}}
\nc{\tr}{\op{tr}}

\nc{\cpt}{{\op{cpt}}}
\nc{\Dol}{{\op{Dol}}}
\nc{\DR}{{\op{DR}}}
\nc{\B}{{\op{B}}}
\nc{\Hod}{{\op{Hod}}}

\nc{\operlimits}[1]{\mathop{\mathchoice{\mbox{\rm #1}}{\mbox{\rm #1}}
{\mbox{\rm \scriptsize #1}}{\mbox{\rm \tiny #1}}}}

\nc{\medoplus}{\operlimits{$\bigoplus$}}
\nc{\Coeff}{\operlimits{Coeff}}


\nc{\al}{\alpha}
\nc{\be}{\beta}
\nc{\ep}{\varepsilon}
\nc{\ga}{\gamma}
\nc{\Ga}{\Gamma}
\nc{\la}{\lambda}
\nc{\La}{\Lambda}
\nc{\si}{\sigma}
\nc{\Sig}{{\Gamma}}
\nc{\Om}{\Omega}
\nc{\om}{\omega}


\nc{\elll}{l}
\nc{\tag}{}

\nc{\down}[1]{{\phantom{\scriptstyle #1}
\hbox{$\left\downarrow\vbox to
    9.5pt{}\right.\nulldelimiterspace=0pt \mathsurround=0pt$}
\raisebox{.4ex}{$\scriptstyle #1$}}}
\nc{\leftdown}[1]{
\hbox{$\raisebox{.4ex}{$\scriptstyle #1$}\left\downarrow\vbox to
    9.5pt{}\right.\nulldelimiterspace=0pt \mathsurround=0pt$}
    \phantom{\scriptstyle #1}}
\nc{\bigdowneq}{{\Big| \! \Big|}}

\nc{\slantbold}[1]{\mathchoice{{\hspace{-.2ex}\mbox{\boldmath
        $#1$}\hspace{-.2ex}}}{{\hspace{-.2ex}\mbox{\boldmath
        $#1$}\hspace{-.2ex}}}{{\hspace{-.2ex}\mbox{\boldmath
        $\scriptstyle
        #1$}\hspace{-.2ex}}}{{\hspace{-.2ex}\mbox{\boldmath
        $\scriptscriptstyle #1$}\hspace{-.2ex}}}  }
\nc{\E}{\slantbold{E}}

\nc{\mc}{\mathcal}
\nc{\ca}{{\mc A}}
\nc{\ce}{{\mc E}}
\nc{\ci}{{I}}
\nc{\cl}{{\mc L}}
\nc{\co}{{\mc O}}
\nc{\cu}{{\mc U}}
\nc{\cv}{{\mc V}}
\nc{\F}{{\mc F}}
\nc{\G}{{\mc G}}
\renewcommand{\L}{{\mc L}}
\nc{\M}{{M}}
\nc{\N}{{\mc N}}
\renewcommand{\O}{{\mc O}}

\renewcommand{\c}{c}

\nc{\GL}[1]{{{\rm GL}_{#1}}}
\nc{\SL}[1]{{{\rm SL}_{#1}}}
\nc{\PGL}[1]{{{\rm PGL}_{#1}}}
\nc{\SU}[1]{{{\rm SU(}#1{\rm )}}}
\nc{\U}[1]{{{\rm U(}#1{\rm )}}}
\nc{\Sp}[1]{{{\rm Sp(}#1{\rm )}}}
\nc{\cx}{{\C^\times}}
\nc{\gli}{{\C^\times}}

\nc{\bino}[2]{\mbox{\Large $#1 \choose #2$}}

\nc{\half}{\mathchoice{{\textstyle \frac{\scriptstyle 1}{\scriptstyle
2}}} {{\textstyle\frac{\scriptstyle 1}{\scriptstyle 2}}}
{\frac{\scriptscriptstyle 1}{\scriptscriptstyle 2}}
{\frac{\scriptscriptstyle 1}{\scriptscriptstyle 2}}}
\nc{\quarter}{\mathchoice{{\textstyle \frac{\scriptstyle
1}{\scriptstyle 4}}} {{\textstyle\frac{\scriptstyle 1}{\scriptstyle
4}}} {\frac{\scriptscriptstyle 1}{\scriptscriptstyle 4}}
{\frac{\scriptscriptstyle 1}{\scriptscriptstyle 4}}}
\nc{\ratio}[2]{\mathchoice{ {\textstyle \frac{\scriptstyle
#1}{\scriptstyle #2}}} {{\textstyle\frac{\scriptstyle #1}{\scriptstyle
#2}}} {\frac{\scriptscriptstyle #1}{\scriptscriptstyle #2}}
{\frac{\scriptscriptstyle #1}{\scriptscriptstyle #2}}}

\nc{\stack}[2]{
{\def\arraystretch{.7}\def\arraycolsep{0pt}\begin{array}{c}
\scriptstyle #1 \\ \scriptstyle #2 \end{array}\def\arraystretch{1}} }


\nc{\Left}[1]{\hbox{$\left#1\vbox to
    10.5pt{}\right.\nulldelimiterspace=0pt \mathsurround=0pt$}}
\nc{\Right}[1]{\hbox{$\left.\vbox to
    10.5pt{}\right#1\nulldelimiterspace=0pt \mathsurround=0pt$}}
\nc{\LEFT}[1]{\hbox{$\left#1\vbox to
    15.5pt{}\right.\nulldelimiterspace=0pt \mathsurround=0pt$}}
\nc{\RIGHT}[1]{\hbox{$\left.\vbox to
    15.5pt{}\right#1\nulldelimiterspace=0pt \mathsurround=0pt$}}

\nc{\lp}{\raisebox{-.1ex}{\rm\large(}\hspace{-.2ex}}
\nc{\rp}{\hspace{-.2ex}\raisebox{-.1ex}{\rm\large)}}

\nc{\emb}{\hookrightarrow}
\nc{\jay}{j}
\nc{\lrow}{\longrightarrow}
\nc{\llrow}{\hbox to 25pt{\rightarrowfill}}
\nc{\lllrow}{\hbox to 29pt{\rightarrowfill}}
\nc{\mbar}{\overline{M}}
\nc{\dbar}{\overline{\partial}}
\nc{\sans}{\,\backslash\,}
\nc{\st}{\, | \,}

\nc{\bee}{{\bf E}}
\nc{\bphi}{{\bf \Phi}}

\nc{\rs}{C}

\nc{\T}{T}
\nc{\Tor}{{\mc T}}
\nc{\Est}{E_{\op{st}}}

\nc{\g}{{\mathfrak g}}
\nc{\gp}[2]{{{\rm #1(}#2{\rm )}}}
\nc{\alg}[2]{{{\mathfrak #1}(#2)}}

\catcode`\@=12



\noindent
{\LARGE \bf Mirror symmetry, Langlands duality, \smallskip \\
and the Hitchin system}
\bigskip \\ 
{\bf Tam\'as Hausel }\\
Department of Mathematics, University of California, 
Berkeley, Calif. 94720\smallskip \\
{\bf Michael Thaddeus } \\
Department of Mathematics, Columbia University, 
New York, N.Y. 10027  
\renewcommand{\thefootnote}{}
\footnotetext{ T.H. supported by a Miller Research Fellowship at the
University of California, Berkeley, 1999--2002.}
\footnotetext{ M.T. supported by NSF grants DMS--98--08529 and DMS--00--99688.}

\medskip

{\small
\begin{quote}
\noindent {\em Abstract.}
Among the major mathematical approaches to mirror symmetry are those
of Batyrev-Borisov and Strominger-Yau-Zaslow (SYZ).  The first is
explicit and amenable to computation but is not clearly related to the
physical motivation; the second is the opposite.  Furthermore, it is
far from obvious that mirror partners in one sense will also be mirror
partners in the other.  This paper concerns a class of examples that
can be shown to satisfy the requirements of SYZ, but whose Hodge
numbers are also equal.  This provides significant evidence in support
of SYZ.  Moreover, the examples are of great interest in their own
right: they are spaces of flat $\SL{r}$-connections on a smooth curve.
The mirror is the corresponding space for the Langlands dual group
$\PGL{r}$.  These examples therefore throw a bridge from mirror
symmetry to the duality theory of Lie groups and, more broadly, to the
geometric Langlands program.

\end{quote}
}

When it emerged in the early 1990s, mirror symmetry was an aspect of
theoretical physics, and specifically a duality between quantum field
theories.  Since then, many people have tried to place it on a mathematical
foundation.  Their labors have built up a fascinating but somewhat
unruly subject.  It describes some sort of relation between pairs of
Calabi-Yau spaces, but there are several quite different formulations
of this relation, with no strong links between them.  Notable among
these are the toric approach of Batyrev-Borisov \cite{b,batboy}, leading
to a very large class of examples whose Hodge numbers behave as
desired, and the symplectic approach of Strominger-Yau-Zaslow \cite{syz},
hereinafter SYZ\@.  The latter is inspired by the original physics,
and holds out the remarkable promise of connecting mirror symmetry to
the theory of integrable systems.  But it is extremely difficult to
find examples.

This paper aims to describe certain pairs of Calabi-Yaus --- namely,
moduli spaces of flat connections on a curve --- which exhibit mirror
symmetry phenomena in two different senses: first, they satisfy the
requirements of SYZ, and second, their Hodge numbers behave more or
less as expected.  As far as we know, these are the first non-trivial
examples of SYZ mirror partners of dimension greater than 2, so they
significantly corroborate the SYZ theory.

Furthermore, our examples relate mirror symmetry to another one of the
great dualities of mathematics: the {\em Langlands duality} on Lie groups.
If $\hat G$ is the Langlands dual of a reductive group $G$, then the
pairs we study are spaces of flat connections on the same base curve
with structure groups $G$ and $\hat G$.  These spaces are basic
objects of study in the geometric Langlands program, which has many
possible points of contact with mirror symmetry.  (For example,
although we do not discuss it here, equivalence of derived categories
of coherent sheaves plays a prominent part in both.)  In the present
paper we confine ourselves throughout to the case $G = \SL{r}$, and
ultimately to the case $G = \SL{2}$ or $\SL{3}$, but we hope and
expect that the mirror relationship holds more generally.  

The original reason for suspecting that our moduli spaces might be
mirror partners was that they comprise dual pairs of {\em
hyperk\"ahler integrable systems}.  The hyperk\"ahler metric and the
collection of Poisson-commuting functions determining the integrable
system were constructed in two seminal papers of Hitchin in the late
1980s \cite{duke,lms}.  These structures automatically produce a
family of special Lagrangian tori on the moduli spaces, which is a key
requirement of SYZ\@.  Moreover, the families of tori on the $\SL{r}$
and $\PGL{r}$ moduli spaces are dual in the appropriate sense, which
is the other requirement of SYZ\@.  The only tricky point is to extend
this story to the moduli spaces of bundles of nonzero degree $d$,
which are technically much easier to deal with when $d$ and $r$ are
coprime.

To deal with this ``twisted'' case, our moduli spaces alone are not
enough: they must be endowed with extra structures, which physicists
call {\em $B$-fields} and mathematicians call {\em flat unitary
gerbes}.  These appear whenever mirror symmetry is formulated in
sufficient generality.  In our case they arise in a particularly
natural way, and indeed they are necessary for things to work properly
when the degree is nonzero.  For instance, as we see in \S\ref{triv},
our case satisfies not the original formulation of SYZ, but rather an
extension proposed by Hitchin \cite{lectures} to Calabi-Yaus with
$B$-fields, of which no examples were previously known.  Likewise, the
Hodge numbers in our case must be evaluated in a generalized sense
involving the $B$-field.  We explain in \S\ref{hodge} how to do this,
adapting the notion of {\em stringy mixed Hodge numbers} as they
appear e.g.\ in Batyrev-Dais \cite{bd}.  These in turn are hybrids of
the {\em stringy Hodge numbers} of Vafa \cite{vafa} and Zaslow
\cite{zas} with the {\em mixed Hodge numbers} of Deligne
\cite{del2,del3}.

Perhaps this is the moment to confess that the relationship between
the Hodge numbers of our mirror partners is not the usual one.  The
familiar identity between Hodge numbers of mirror partners $M$ and
$\hat M$ is of the form $h^{p,q}(\hat M) = h^{\dim M - p, q}(M)$.  We
will see, however, that our mirror partners satisfy an identity of a
simpler form: just $h^{p,q}(\hat M) = h^{p,q}(M)$.  This seems to
reflect the fact that they are hyperk\"ahler and noncompact.  At any
rate, compact hyperk\"ahler manifolds (and orbifolds) satisfy
$h^{p,q}(M) = h^{\dim M - p, q}(M)$, and hence we expect $h^{p,q}(\hat
M) = h^{p,q}(M)$ for compact hyperk\"ahler mirror partners.
Apparently this relationship persists in the noncompact case, even
though the familiar mirror identity does not.

A physical explanation of this based on the original quantum field
theory would be gratifying.  But we must also bear in mind that, for
noncompact varieties, the mixed Hodge numbers, and hence our Hodge
numbers, depend on the algebraic structure.  (Indeed, the spaces of
representations of the fundamental group --- what Simpson \cite{simp}
calls the Betti spaces --- are analytically but not algebraically
isomorphic to our spaces, and their Hodge numbers will in general be
different.)  This seems hard to explain from a physical point of view.
It might be preferable to work with some notion of Hodge numbers
depending on the metric and not the algebraic structure.

Nevertheless, the equality of Hodge numbers that we uncover is
striking and totally unexpected from a mathematical viewpoint.  At any
rate, it follows from the equality of terms contributed by loci in the
moduli space which seem to be completely unrelated to one another.
They are fixed loci of natural group actions, but on one side, the
group acting is $\cx$, while on the other it is a finite abelian group
$\Ga$.  So our result illustrates both the power and the mystery of
mirror symmetry.

Here is a sharper outline of the contents of the paper.  The
first two sections review the known facts we will need: \S\ref{syz}
covers Calabi-Yaus, gerbes, and the proposal of SYZ, while
\S\ref{higgs} covers Higgs bundles, flat connections, and the Hitchin
system.  The next section is devoted to the proof of our first main
result, Theorem \re{ll}, showing that the moduli spaces of flat
connections on a curve with structure groups $\SL{r}$ and $\PGL{r}$
are SYZ mirror partners.

The rest of the paper is devoted to the evaluation of Hodge numbers
for these spaces.  In \S\ref{hodge} we define the appropriate notion
of Hodge numbers: stringy mixed Hodge numbers with coefficient system
provided by a flat unitary gerbe!  This enables us to state our main
conjecture, Conjecture \re{jj}, on the equality of the Hodge numbers
for the $\SL{r}$ and $\PGL{r}$ spaces, which we then proceed to prove for
$r=2$ and $3$.

It is much easier to work with Higgs bundles than flat connections,
because of the algebraic $\cx$-action on the moduli space.  So we
begin our proof by showing in \S\ref{equal} that these two moduli
spaces have the same Hodge numbers, and thereafter we work exclusively
with the space of Higgs bundles.  In \S\ref{I} we describe (following
Narasimhan-Ramanan \cite{nr}) the fixed points of the action on the
$\SL{r}$ moduli space of the group $\Ga$ of $r$-torsion points in the
Jacobian, and in \S\ref{II} we use this to compute the Hodge numbers
of the $\PGL{r}$ space.  Then in \S\ref{IV} we describe (following
Hitchin \cite{lms} and Gothen \cite{goth}) the fixed points of the
action of $\cx$, and in \S\ref{V} we use this to compute the Hodge
numbers of the $\SL{r}$ space --- in sufficient detail that, for $r=2$
and $3$, we get a complete answer.

The main results of this paper were announced in a note in 2001
\cite{ht}.  The Proposition and Theorem 3 in that announcement
correspond roughly to Theorems \re{ll} and \re{mm} in the present
work.  But the latter results actually represent substantial
improvements: in particular, the meaning and function of the $B$-field
have been greatly clarified.  For example, Theorem 3 of the
announcement refers to stringy Hodge numbers with discrete torsion;
although the numbers turn out to be the same, we now understand that
the canonical $B$-field, as defined in \S\ref{triv}, is a flat gerbe
which may not come from discrete torsion.  Theorems 1 and 2 of the
announcement concern the spaces of flat connections on {\em
punctured\/} curves, or equivalently, of {\em parabolic\/} Higgs
bundles; once the $B$-field is properly understood, this is mostly
parallel to the present case, and we intend to treat it elsewhere.

One word about terminology: we use torsors liberally in the paper,
both for sheaves of groups and for group schemes, so here is a
definition.  A {\em torsor} for a sheaf of groups $T$ over a base $X$
is a sheaf of $T$-spaces over $X$ which is locally isomorphic to $T$
as a sheaf with $T$-action.  The same definition holds if sheaves are
replaced by schemes, or even by families in a $C^\infty$ sense.  
If $X$ is a point, then a $T$-torsor is a {\em principal homogeneous
space} for $T$: it can be identified with $T$ up to the choice of a
basepoint.  In this paper, the relevant groups are always abelian.

{\em Acknowledgements.}  We are very grateful to Nigel Hitchin, who
suggested the germ of the idea for this work as long ago as 1996; and
to Pierre Deligne, Ron Donagi, Dennis Gaitsgory, Tony Pantev, and
Bal\'azs Szendr\H oi for helpful remarks.  We also thank Cumrun Vafa
for drawing our attention to an earlier paper \cite{bjsv} treating
related ideas in a physical context.

\bit{Strominger-Yau-Zaslow} 
\label{syz}

\noindent {\bf Calabi-Yau manifolds and B-fields.}
We take a {\em Calabi-Yau} manifold to be a complex manifold equipped
with a Ricci-flat K\"ahler metric.  On a Calabi-Yau manifold of
complex dimension $n$, parallel transport defines on any simply
connected neighborhood a covariant constant holomorphic $n$-form
$\Om$, unique up to a scalar.  Usually --- as for example when $M$
itself is simply connected --- this form is defined globally, and we
assume this for simplicity.

Mirror symmetry is supposed to relate two such Calabi-Yau manifolds
$M$ and $\hat M$, interchanging the deformation spaces of the K\"ahler
and complex structures.  However, the K\"ahler forms are real 2-forms
of type $(1,1)$; to allow the K\"ahler deformations to be complex, we
choose auxiliary fields, say $B$ on $M$ and $\hat B$ on $\hat M$,
which are in some sense imaginary parts for the K\"ahler forms.
Exactly where the $B$-field takes values is not entirely clear in the
physics literature, but following a suggestion of Hitchin
\cite{lectures} we will take it to be an element of $H^2(M,\U{1})$, or
an isomorphism class of flat unitary gerbes.  By this we mean the
following.

\noindent {\bf Gerbes and their trivializations.}
Let $T$ be a sheaf of abelian groups over a variety $M$ (with the
complex or \'etale topology).  A {\em Picard category\/} is a tensor
category where all objects and all morphisms are invertible.  The
category of $T$-torsors constitutes a sheaf of Picard categories (or
stack) over $M$.  {\em Sheaf of categories\/} here means roughly what
one would expect, but the precise definition is somewhat technical; a
convenient reference is Donagi-Gaitsgory \cite{dg}.

A $T$-{\em gerbe} is a sheaf of categories which is a torsor over this
sheaf.  That is, the sheaf consisting of $T$-torsors must act on the
gerbe, and be locally equivalent to it as a sheaf with this action.
For us, $T$ will usually be the sheaf of locally constant functions
with values in $\U{1}$; then $\U{1}$-torsors are flat unitary line
bundles, and we refer to $\U{1}$-gerbes as {\em flat unitary gerbes}.

An {\em isomorphism} of $T$-gerbes is an equivalence of sheaves of
categories as torsors over the sheaf of $T$-torsors.  An {\em
automorphism} is a self-isomorphism; since a gerbe is a torsor over
the sheaf of $T$-torsors, its automorphisms are identified with
sections of that sheaf, that is, with $T$-torsors, acting by
tensorization.  A {\em trivialization} of a $T$-gerbe is an
isomorphism to the trivial gerbe.  Two trivializations $z,z'$ are {\em
equivalent} if the automorphism $z' \circ z^{-1}$ is given by
tensorization with a trivial $T$-torsor.  The space of equivalence
classes of trivializations of a (trivial) gerbe $B$ will be denoted
$\Triv^T(M,B)$; it is naturally an $H^1(M,T)$-torsor over a point.

The key result on gerbes is due to Giraud \cite{bryl,giraud}: it
asserts that isomorphism classes of $T$-gerbes are in one-to-one
correspondence with $H^2(M, T)$.  Indeed, to construct a \v Cech
cocycle $x \in C^2(M, T)$ from a gerbe, choose a cover $\{ U_\al \}$
such that the gerbe is trivialized on each $U_\al$.  The overlaps are
then given by tensorizations by $T$-torsors $L_{\al,\be}$, with
$L_{\al,\be} \otimes L_{\be,\ga} \otimes L_{\ga,\al}$ canonically
trivialized on the triple overlaps.  After refining the cover if
necessary, trivialize each $L_{\al,\be}$, and then compare with the
canonical trivializations on the triple overlaps to get the cocycle.
In this setting, a trivialization can be regarded as a cochain $y \in
C^1(M, T)$ with $dy=x$, and two trivializations are equivalent if they
differ by an exact cocycle.  Then it is clear why equivalence classes
of trivializations form an $H^1(M,T)$-torsor.

\noindent {\bf Orbifolds.}  Strictly speaking, the mirror of a
Calabi-Yau manifold may not be a manifold, but rather an orbifold.
The notion of a Calabi-Yau orbifold is defined in Appendix A of
Cox-Katz \cite{cz}, and on such orbifolds, gerbes may be defined much as
line bundles are.

For the present purposes no difficult theory is needed, as the
orbifolds we encounter are all global quotients of Calabi-Yau
manifolds by the actions of finite groups.  If $M = X/\Ga$ is a
quotient of this kind, and $T$ is a sheaf over $X$ to which the action
of $\Ga$ lifts, then a $T$-gerbe on $M$ is simply a $T$-gerbe on $X$
equipped with a lifting of the $\Ga$-action.

For example, if $X$ is a point, then a $\Ga$-equivariant $\U{1}$-gerbe
is a homomorphism from $\Ga$ to the category of $\U{1}$-torsors over a
point, which is nothing but a central extension of $\Ga$ by $\U{1}$.
Such extensions are classified up to isomorphism by the group
cohomology $H^2(\Ga,\U{1})$.  In the physics literature, this last
group is called the {\em discrete torsion\/} of $\Ga$ \cite{vw}; in the
mathematics literature, it is called the {\em Schur multiplier} \cite{karp}.

\noindent {\bf Strominger-Yau-Zaslow.}  
With all this understood, the proposal of Strominger-Yau-Zaslow can be
described as follows.  

A torus $L$ of real dimension $n$ embedded in a Calabi-Yau
$n$-orbifold is said to be {\em special Lagrangian} if $\om|_L =0$ and
$\Im \Om |_L =0$.

Two Calabi-Yau $n$-orbifolds $M$ and $\hat M$,
equipped with flat unitary gerbes $B$ and $\hat B$, respectively, 
are said to be {\em SYZ mirror partners\/} if there exist an
orbifold $N$ of real dimension $n$ 
and smooth surjections $\mu:M \to N$ and $\hat \mu:
\hat M\to N$ such that for every $x \in N$ 
which is a regular value of $\mu$ and $\hat \mu$,
the fibers $L_x = \mu^{-1}(x) \subset M$ 
and $\hat L_x = \hat \mu^{-1}(x) \subset \hat M$ 
are special Lagrangian tori which are dual to each other in the
sense that there are smooth identifications
$$L_x = \Triv^\U{1} (\hat L_x, \hat B)$$ and 
$$\hat L_x = \Triv^\U{1} (L_x, B)$$
depending smoothly on $x$.
The requirement that the identifications only be smooth is rather
weak, but it is unclear what a stronger condition ought to be.
Certainly isometry is too strong.  

\noindent {\bf The hyperk\"ahler case.}  
Constructing special Lagrangian tori is usually very difficult.  But
suppose that $M$ is a {\em hyperk\"ahler} manifold: that is, it has a
metric which is simultaneously K\"ahler with respect to three complex
structures $J_1,J_2,J_3$ satisfying the commutation relations of the
quaternions $i,j,k$.  
Let $\om_1,\om_2, \om_3$ be the corresponding K\"ahler forms.  Then
$\om_2+ i \om_3$ is a complex symplectic form on $M$ which is
holomorphic with respect to $J_1$.  The associated volume form gives a
covariant constant trivialization of the canonical bundle, which shows
that $M$ is Ricci-flat and hence Calabi-Yau with respect to $J_1$, and
by permuting the indices, with respect to all three complex
structures.  

In this case, it is easy to see that any complex submanifold $L\subset
M$ which is complex Lagrangian with respect to $J_1$ is special
Lagrangian with respect to $J_2$ \cite{pisa}.  So the desired family
of special Lagrangian tori can be found by holomorphic methods: first
find holomorphic maps $\mu$ and $\hat \mu$ whose generic fibers are
complex Lagrangian tori, then perform a {\em hyperk\"ahler rotation},
that is, change to a different complex structure.  The Hitchin system,
to be described below, gives holomorphic maps of precisely this kind
on a hyperk\"ahler manifold.

\bit{Higgs bundles and local systems}
\label{higgs}

\noindent {\bf Review of the basic facts.}
Let us recall some of the theory of Higgs bundles and local systems on
curves, as developed by Hitchin \cite{lms} and Simpson \cite{simp}.

Let $C$ be a smooth complex projective curve of genus $g$.  It will be
convenient to fix a basepoint $c \in C$.  A {\em Higgs bundle} is a
pair $(E,\phi)$ consisting of a vector bundle $E$ over $C$ and a
section $\phi \in H^0(C,\End E \otimes K)$, where $K$ is the canonical
bundle.  It is {\em stable} 
if all proper subbundles $F \subset E$ with $\phi(F) \subset F
\otimes K$ satisfy 
$\deg F / \rk F < 
\deg E / \rk E$.  
Hitchin and Simpson then prove the following.
(The subscripts Dol, DR and Hod are Simpson's notation; they honor
Dolbeault, De Rham, and Hodge respectively.)

\begin{itemize}
  
\item There exists a smooth, quasi-projective moduli space $M^d_\Dol$
  of stable Higgs bundles of rank $r$ and degree $d$.  
  
\item There exists a smooth, quasi-projective moduli space $M^d_\DR$
  of irreducible local systems (that is, flat vector bundles) on $C
  \sans \{ c\}$ of rank $r$, with holonomy $e^{2\pi id/r}$ around $c$.
  
\item These two spaces are naturally diffeomorphic; indeed, there
  exists an isosingular family $M^d_\Hod$ over the affine line whose
  zero fiber is $M^d_\Dol$, but whose fiber over every other closed
  point is $M^d_\DR$.
  
\item There is a $\cx$-action on $M^d_\Hod$ lifting the standard
  action on the affine line, and restricting to $M^d_\Dol$ as $t
  \cdot (E,\phi) = (E, t \phi)$.
  
\item There is a Riemannian metric on $M^d_\Dol$ for which the
  Dolbeault and De Rham complex structures form part of a
  hyperk\"ahler structure.

\end{itemize}

\smallskip

\noindent {\bf More general structure groups.}
If vector bundles are replaced by principal bundles, the whole theory
generalizes without difficulty.  Higgs bundles and local systems make
sense, their moduli make sense, and even the spaces $M_\Hod$ make
sense.  For example, the right notion of a {\em principal Higgs
bundle\/} consists of a principal bundle $E$ and a section $\phi \in
H^0(C,\ad E \otimes K)$.  Simpson explains why there exists a moduli
space of principal Higgs bundles, stable in the appropriate sense.
However, for the purposes of this paper we only need structure groups
$\GL{r}$, $\SL{r}$, and $\PGL{r}$, so we make do with the direct
descriptions of the moduli spaces below.  It is easy to check that
these descriptions agree with Simpson's definitions, but it is even
easier to regard these descriptions as definitions themselves.  In
each case, we describe $M_\Hod$; $M_\Dol$ and $M_\DR$ are the zero and
nonzero fibers, respectively.
\begin{itemize}

\item Let $M^d_\Hod(\GL{r}) = M^d_\Hod$ as defined above.  

\item Let $M^d_\Hod(\SL{r})$ be the inverse image of a smooth
algebraic section $s \subset M^d_\Hod(\cx)$ under the morphism $\det:
M^d_\Hod(\GL{r}) \to M^d_\Hod(\cx)$ induced by the determinant
representation of $\GL{r}$.  It is convenient to take $s(0) =
(\O(dc),0) \in M^d_\Dol(\cx)$, where $c \in C$ is the basepoint; then
$M^d_\Dol(\SL{r})$ parametrizes stable Higgs bundles $(E,\phi)$ with
$\La^r E \cong \O(dc)$ and $\tr \phi = 0$.

\item Let $M^d_\Hod(\PGL{r})$ be the geometric quotient of
$M^d_\Hod(\SL{r})$ by the group scheme $\Sig = \Pic^0 C[r]$ consisting
of isomorphism classes of line bundles whose $r$th power is trivial,
acting by tensorization.
 
\end{itemize}

In the second item, the existence of a smooth section $s$ follows, for
example, from the Bia\l ynicki-Birula decomposition theorem \cite{bb,
flip}, taking a $\cx$-orbit whose closure contains $(\O(dc),0)$.  To
see that every section $s$ gives the same space $M^d_\Hod(\SL{r})$ up
to isomorphism, notice that tensorization makes $M^d_\Hod(\cx)$ into a
$M^0_\Hod(\cx)$-torsor over the affine line, so any two sections
differ by multiplication by a section of the family of groups
$M^0_\Hod(\cx)$.  Since the $r$th power map is \'etale on
$M^0_\Hod(\cx)$, and the base is simply connected, there exists a
smooth $r$th root of this section, which can be used to tensorize
objects in $M^d_\Hod(\GL{r})$.

Incidentally, the use of the notation $M^d_\Hod(\SL{r})$ is perhaps
slightly misleading, since the objects it parametrizes do not have
structure group $\SL{r}$ unless $d=0$ and $s$ is the trivial section.

Everything asserted before for Higgs bundles remains true in this more
general setting, except that the $\PGL{r}$ moduli spaces are not
smooth; rather, they are hyperk\"ahler orbifolds.

\smallskip

\noindent {\bf The Hitchin system.}
On each of the Dolbeault spaces $M^d_\Dol(G)$, there exists a
completely integrable Hamiltonian system, the so-called {\em Hitchin
system} or {\em Hitchin map}.  It is a morphism $\mu$ from
$M^d_\Dol(G)$ to an affine space $V_G$ of half the dimension.  Here
$V_{\GL{r}} = \bigoplus_{i=1}^n H^0(C,K^i)$ and $V_{\SL{r}}
=V_{\PGL{r}} = \bigoplus_{i=2}^n H^0(C,K^i)$.  The morphism is
evaluated on a Higgs bundle $(E,\phi)$ by applying to $\phi$ the
invariant polynomials on the Lie algebra $\mathfrak g$.  Hitchin shows
that $\mu$ is proper when $r$ and $d$ are coprime.  He also shows that
the fiber over a general point is complex Lagrangian and is (a torsor
for) an abelian variety.

This is exactly the situation discussed at the end of \S\ref{syz}.
Consequently, for any integers $d,e \in \Z$, the De Rham spaces
$M^d_\DR(\SL{r})$ and $M^e_\DR(\PGL{r})$ carry families of special
Lagrangian tori over the same base, just as the SYZ definition
requires.  All that remains to be verified is the statement about
duality of the tori.  We will establish this in the next section, but
first we need to review Hitchin's description of the fibers of
$\mu$ in more detail.

An element of $V_G$ is given by sections $\beta_i \in H^0(C,K^i)$ for
$i = 1$ to $n$ (taking $\beta_1 = 0$ in the $\SL{r}$ and $\PGL{r}$
cases).  The equation
\beq
\label{nn}
z^n + \beta_1 z^{n-1} + \beta_2 z^{n-2} + \cdots + \beta_n = 0,
\eeq
where $z$ lies in the total space of $K$, defines a curve $\pi: \tilde
C \to C$ inside this total space, called the {\em spectral cover}.
For $(\beta_i)$ in the Zariski open set $U \subset V_G$ where $\tilde
C$ is smooth, $\mu^{-1}(\beta_i)$ can be canonically identified as
follows \cite{lms}. 
\begin{itemize}
  
\item When $G = \GL{r}$, it is $\tilde J^d = \Pic^d \tilde C$.  This
  can be regarded as (the fiber of) a $\tilde J^0$-torsor over $U$.
  
\item When $G = \SL{r}$, it is $P^d = \Nm^{-1}(\O(dc))$, the {\em
  generalized Prym variety}.  Here $\Nm: \Pic^d \tilde C \to \Pic^d
  C$ is the norm map (see e.g.\ Arbarello et al.\ \cite[App.\ B]{acgh}
  for a definition and basic properties).  This is a $P^0$-torsor over $U$.

\item When $G = \PGL{r}$, it is $\hat P^d = P^d/\Sig$, the quotient of
  the Prym by the action of $\Sig$ on $\M^d_\Dol(\SL{r})$, which
  preserves it.  This is a $\hat P^0$-torsor over $U$.

\end{itemize}

The next two lemmas explain how these torsors are related to one another. 

\bs{Lemma}
\label{qq}
Let $J^0 = \Pic^0 C$.  Then there is a natural isomorphism
$$\tilde J^d \cong \frac{P^d \times J^0}{\Ga}$$
under which $\Nm$ corresponds to the projection to $J^0/\Ga$ followed
by the isomorphism $J^0/\Ga \to J^0$ given by taking $-r$th powers.
\es

\pf.  Certainly there is a morphism $P^d \times J^0 \to \tilde J^d$
given by $(L,M) \mapsto L \otimes \pi^* M^{-1}$, whose composition
with $\Nm$ is $(L,M) \mapsto M^{-r}$.  This morphism is invariant
under the action of $\Ga$ by tensorization on both factors.  It
therefore suffices to show that $P^0 \cap \pi^*J^0 = \pi^* \Ga$ and
that $\pi^*: J^0 \to \tilde J^0$ is injective.

Since $\Nm \pi^* L = L^r$, certainly $\ker \pi^* \subset \Ga$ and
$P^0 \cap \pi^* J^0 = \ker \Nm \cap \pi^* J^0 = \pi^* \Ga$,
which proves the first assertion.  

For the second, suppose $L \in \Ga$ has order $k$ and satisfies $\pi^*
L \cong \O$.  The isomorphism $L^k \cong \O$ determines a $k$-valued
multisection of $L$; regard this as a cover of $C$.  Then $\pi:
\tilde C \to C$ must factor through this cover: indeed, the
trivialization of $\pi^* L$ gives a trivialization of $\pi^* L^k$, so
(after multiplication by an overall constant) it lies in the pullback
of, and so defines a morphism to, the multisection.

However, the only connected unbranched cover of $C$ through which
$\pi$ factors is the trivial cover.  This is clear when $\pi$ has a
point of total ramification, which occurs when all the coefficients
$\beta_i$ of \re{nn} have a common zero.  But it is also clearly
invariant under continuous deformations, hence true everywhere on the
connected base $U \subset V_G$.

Therefore $k=1$, so $\pi^*$ has trivial kernel.  \fp

\bs{Lemma}
\label{oo}
The dual of $P^0$ is $\hat P^0 = P^0 / \Ga$.
\es

\pf.  Just dualize the short exact sequence 
$$ 0 \lrow P^0 \lrow \tilde J^0 \stackrel{\Nm}{\lrow} J^0 \lrow 0$$
to get 
$$ 0 \lrow J^0 \stackrel{\pi^*}{\lrow} \tilde J^0 \lrow P^0 \lrow 0,$$
where $\hat P^0 = \tilde J^0 /J^0 = P^0 /\Ga$ by the previous lemma.  \fp

\bit{Trivializations of the \boldmath $B$-fields}
\label{triv}

With the prerequisites complete, we proceed in this section to our
first goal.  This is to show that, when equipped with certain
$B$-fields, the De Rham moduli spaces with structure groups $G$ and
$\hat G$ are SYZ mirror partners.  Our expectation is that this will hold
true for any reductive $G$, but at present we confine ourselves to the
case $G = \SL{r}$, $\hat G = \PGL{r}$.

In fact, we deduce the smooth identification of special Lagrangian
tori on the De Rham spaces, called for by SYZ, from a stronger
statement: a holomorphic identification of complex Lagrangian tori on
the Dolbeault spaces.  The two are related by hyperk\"ahler rotation
as discussed at the end of \S\ref{syz}.  Moreover, since the smooth
parts of these two spaces are diffeomorphic, flat unitary gerbes on
them can be identified.  Therefore, in this section, we work
exclusively with the Dolbeault space and, for brevity, denote the
stable part of $\M^d_\Dol(\SL{r})$ simply by $\M^d_\Dol$.

We will work over $U$, the Zariski open set in the range of the
Hitchin map $\mu$ where the spectral cover $\tilde C$ is smooth.  The
four torsors over $U$ that concern us are $\tilde J^d$, $P^d$, and
$\hat P^d$, as defined in \S\ref{higgs}, plus $J^d = \Pic^d C$, which
we regard as a trivial $J^0$-torsor over $U$.

Any of the methods used to construct universal families on the moduli
space of ordinary stable bundles adapt without change to the space of
Higgs bundles $\M^d_\Dol$.  For example, one could use descent and the
geometric invariant theory construction of $\M^d_\Dol$ given by
Nitsure \cite{nitsure}.  Provided that $r$ and $d$ are coprime, one
gets a bona fide universal Higgs pair $(\bee, \bphi) \to \M^d_\Dol
\times C$.  However, as for stable bundles, the scalars, acting as
automorphisms, provide an obstruction to the existence of $\bee$ when
$r$ and $d$ are not coprime.  The best we can do in general is to
construct a universal projective bundle and a universal endomorphism
bundle, abusively denoted $\P \bee$ and $\End \bee$ even though $\bee$
does not exist, and a universal Higgs field $\bphi \in H^0(\End \bee
\otimes K)$.  Then the restriction $\P \bee|_{\M^d_\Dol \times \{ \c
  \} }$ to the basepoint in $C$ is a projective bundle $\Psi$ on
$\M^d_\Dol$.

Let $B$ be the {\em gerbe of liftings of $\Psi$}, meaning the sheaf of
categories on $M$ taking an \'etale neighborhood to the category of
liftings on that neighborhood of $\Psi$ to an $\SL{r}$-bundle.  Since
any two liftings differ by tensorization by a $\Z_r$-torsor, $B$ is a
$\Z_r$-gerbe.  

\bs{Lemma}
\label{cc}
The restriction of $B$ to each regular fiber $P^d$ of
the Hitchin map is trivial as a $\Z_r$-gerbe.
\es

\pf.  It suffices to show that $\Psi|_{P^d}$ does in fact lift to an
$\SL{r}$-bundle, which can be regarded as a vector bundle with trivial
determinant.  

For any universal bundle $\tilde \L \to P^d \times \tilde C$, the
push-forward $\pi_* \tilde \L \to P^d \times \tilde C$ admits a family
of Higgs fields inducing the inclusion $P^d \subset M$.  Indeed, this
is how one shows that $P^d$ is the fiber of the Hitchin map: see
Hitchin \cite{lms} for details.  So over $P^d \times \{ c \}$ there is
an isomorphism $\P \pi_* \tilde \L \cong \Psi$.

The universal bundle can be normalized so that $\tilde \L|_{P^d \times
\{y \}} \in \Pic^0(P^d)$ for one (hence all) $y \in \tilde C$.  The
determinant of $\pi_* \tilde \L$ over $P^d \times \{c\}$ is isomorphic
to $\bigotimes_{y \in \pi^{-1}(c)} \tilde \L|_{P^d \times \{ y \}}$,
where ramification points are counted with the appropriate
multiplicity.  This has an $r$th root, tensoring by whose inverse will
further adjust the normalization of $\tilde \L$ so that $\det \pi_*
\tilde \L |_{P^d \times \{c \}} \cong \O$, making $\pi_* \tilde
\L|_{P^d \times \{c \}}$ an $\SL{r}$-bundle as desired.  \fp

Now that we know that $B$ restricts trivially to each fiber as a
$\Z_r$-gerbe, and hence as a $\U{1}$-gerbe, it makes sense to examine
the equivalence classes of $\U{1}$-trivializations.  From the discussion of
gerbes in \S\ref{syz}, we know that these form a torsor (in the smooth
category) for $H^1(P^d,\U{1}) \cong \Pic^0 P^d \cong \Pic^0 P^0$, and
from Lemma \re{oo} the latter is $\hat P^0$.

\bs{Proposition}
\label{ee}
For any $d,e \in \Z$, there is a smooth isomorphism of 
$\hat P^0$-torsors $$\Triv^{\U{1}}(P^d,B^e) \cong \hat P^e.$$
\es

\pf.  The isomorphism classes of torsors over a fixed abelian group scheme
themselves form an abelian group in a natural way, and it is easy to
see that $\hat P^e \cong (\hat P^1)^e$ and $\Triv(P, B^e) \cong
(\Triv(P, B))^e$, where the right-hand sides are $e$th powers under
this group operation.  Hence it suffices to take $e=1$.

As seen in the proof of Lemma \re{cc}, the triviality of $B$ on $P^d$
follows from the existence of a universal bundle $\tilde \L \to P^d
\times \tilde C$ with $\det \pi_* \tilde \L$ trivial on $P^d \times \{
c \}$.  Consider the set of isomorphism classes of all such $\tilde
\L$: this parametrizes the equivalence classes of trivializations of
$B$ as a $\Z_r$-gerbe, which is to say, it forms the torsor
$\Triv^{\Z_r} (P^d,B)$.  It is a $\hat P^0[r]$-torsor over $U$, where
$\hat P^0[r] = H^1(P^d,\Z_r)$ are the torsion points of order $r$ in
$\hat P^0$.  This makes sense, since for $L \in \hat P^0 =
\Pic^0(P^d)$, the push-pull formula says
$$\det \pi_*(\pi^* L \otimes \tilde \L) = L^r \otimes \det \pi_* \tilde \L.$$

We are really interested not only in $\Z_r$-trivializations but in all
$\U{1}$-trivializations.  These comprise a torsor for $\hat P^0 =
H^1(P^d,\U{1})$ containing the $\hat P^0[r]$-torsor above, and indeed
this property determines the larger torsor, since it can be identified
with the quotient
$$\frac{\Triv^{\Z_r} (P^d,B) \times \hat P^0}{\hat P^0[r]}.$$

An obvious torsor with this property consists of all universal bundles
$\tilde \L \to P^d \times \tilde C$ with $\tilde \L|_{P^d \times \{ y \} } \in
\Pic^0(P^d)$ for any $y \in \tilde C$.  It therefore suffices to show
that this torsor is isomorphic to $\hat P^1$.

In fact, $\hat P^1 = \tilde J^1 / J^0$, while the torsor of the
previous paragraph is also a quotient by $J^0$, of the torsor
consisting of universal bundles as stated there, except with $P^d$
replaced by $\tilde J^d$.  (The $J^0$-action comes from tensoring by
$\pi^* \Pic J^0 \cong J^0$.)  So it actually suffices to show that the
latter torsor is isomorphic to $\tilde J^1$ as a $\tilde J^0$-torsor.
To do this, we will exhibit morphisms $f_1$ and $f_2$ from $\tilde C$
to the two torsors such that, for any $y,y' \in \tilde C$, $f_1(y') -
f_1(y) = f_2(y') - f_2(y) \in \hat P^0$.  The isomorphism of the two
torsors defined by identifying $f_1(y)$ with $f_2(y)$ is then
independent of $y$, and hence well-defined.

The morphism $f_1$ is simply the Abel-Jacobi map $\tilde C \to \tilde
J^1$.  As for $f_2$, it takes $y$ to the universal bundle whose
restriction to $P^d \times \{y\}$ is trivial.  The equality $f_1(y') -
f_1(y) = f_2(y') - f_2(y)$ then means that the restriction to $y$ of
the universal bundle normalized at $y'$ is the line bundle on $\tilde
J^d$ corresponding to $f_1(y') - f_1(y) \in \tilde J^0 = \Pic^0 \tilde
J^d$.  This follows readily from two well-known facts.  First, that
this universal bundle is of the form $p_2^* L_0\otimes F^* {\mathcal
P}$, where $p_2$ is projection on the second factor, $L_0 \in
\tilde J^d$ is fixed, $\mathcal P$ is the Poincar\'e line bundle, and
$$F: \tilde J^d \times \tilde C \to \tilde J^0 \times \tilde J^0$$
is given by $F(L,y) = \Left(L\otimes L_0^{-1},f_1(y) -
f_1(y')\Right)$.  Second, that the involution of $\tilde J^0 \times
\tilde J^0$ exchanging the two factors takes the Poincar\'e bundle to
its inverse.  \fp

Now turn to the reverse direction.  We need a gerbe $\hat B$ on the
orbifold $\hat \M^d_\Dol = \M^d_\Dol / \Ga$, or equivalently, a
$\Ga$-equivariant gerbe on $\M^d_\Dol$.  This will just be $B$
equipped with a $\Ga$-action, which we define as follows.  Let $L_\ga$
denote the line bundle over $C$ corresponding to $\ga \in \Ga$.  (Here
and throughout, it will prove convenient to distinguish between the
abstract group element $\ga$ and the line bundle $L_\ga$.)  Then $\ga$
acts on $\M^d_\Dol$ by $(E,\phi) \mapsto (E \otimes L_\ga, \phi)$.
This lifts to $\P \bee$: think of $\P \bee$ as the moduli space
parametrizing 1-dimensional subspaces of a stable Higgs bundle, and
observe that tensoring by $L_\ga$ induces a natural transformation.
Hence $\Ga$ acts on $\P \bee$ and on its restriction $\Psi$ to
$M^d_\Dol \times \{ c \}$.  This determines a $\Ga$-action on $B$, the
sheaf of liftings to $\SL{r}$-bundles, making it an equivariant flat
gerbe $\hat B$.

To prove an analogue of Lemma \re{cc} for $\hat B$, we first need a 
technical fact.  Let 
$$ \tilde \Ga = \bigsqcup_{\ga \in \Ga} L_\ga \sans 0$$
be the disjoint union of the total spaces of the line bundles $L_\ga$,
minus their zero sections.

\bs{Lemma}
\label{pp}
This has the structure of a group scheme over $C$ whose fiber at $y
\in C$ is an abelian extension
$$1 \lrow \C^\times \lrow \tilde \Ga_y \lrow \Ga \lrow 0.$$ 
If $\L \to J^0 \times C$ is the universal bundle which is trivial on
$J^0 \times \{ \c \}$, then there is an action over $C$ of $\tilde
\Ga$ on the total space of $\L$, lifting the action of $\Ga$ on $J^0$
by translation, so that the scalars $\C^\times$ act with weight $1$ on
the fibers.  \es

Of course the above extension always splits, but not canonically
except at the basepoint $y=\c$.

\pf.  Let $A$ be an abelian variety (which we will shortly take to be
$J^0$), let $\hat A$ be its dual, and let ${\mathcal P} \to A \times
\hat A$ be the Poincar\'e bundle trivialized on $(0 \times \hat A)
\cup (A \times \hat 0)$, where $0 \in A$, $\hat 0 \in\hat A$ are the
basepoints.  It is well-known, cf.\ Serre \cite[VII 3.16]{serre}, that
$\Ext^1(A,\C^\times) = \hat A$; indeed, $\hat A$ parametrizes a family
of abelian central extensions of $A$ by $\C^\times$.  If this is
regarded as a group scheme over $\hat A$, then its total space is
${\mathcal P} \sans 0$, the complement of the zero section in
$\mathcal P$.  The group operation over $\hat A$ is given by an
isomorphism over $A \times A \times \hat A$ 
\beq
\label{p} 
p_{13}^* {\mathcal P} \otimes p_{23}^* {\mathcal P} 
\cong (m \times 1)^* {\mathcal P}, 
\eeq 
where $p_{13}$ and $p_{23}$ are projections on the relevant factors
and $m: A \times A \to A$ is addition.  This is provided by the
theorem of the cube \cite{mum}.  It can be chosen so that over $0
\times \hat 0$ it is 1.  Associativity requires the commutativity of a
certain diagram of isomorphisms of line bundles on $A \times A \times
A \times \hat A$, but this is automatic since the base is projective
and connected, and the desired commutativity holds automatically over
the base points.

In the same way, the action of the group scheme $A \times \hat A \to
\hat A$ on itself by translating the first factor lifts to an action
of ${\mathcal P} \sans 0 \to \hat A$ on the Poincar\'e bundle
${\mathcal P} \to A \times \hat A \to \hat A$.  Indeed, the action
$({\mathcal P} \sans 0) \times_{\hat A} {\mathcal P} \to {\mathcal P}$
is again given by the isomorphism \re{p}.  The condition that an
action must satisfy is automatic for the same reason as before.

Now return to our curve $C$ with basepoint $\c$, and use the
Abel-Jacobi map to embed it in its Jacobian $J^0$ so that $\c$ maps to
the basepoint $0$.  Let $\tilde \Ga$ be the inverse image of $\Ga
\times C$ in the projection ${\mathcal P} \sans 0 \to J^0 \times J^0$.
Then $\tilde \Ga$ clearly satisfies the desired properties.
The restriction of the Poincar\'e bundle on $J^0 \times J^0$ to $J^0
\times C$ is the universal bundle $\L$; this therefore carries the
desired action.  \fp

\bs{Lemma}
The restriction of $\hat B$ to each regular fiber $\hat P^d$ of the
Hitchin map is trivial as a $\Z_r$-gerbe.
\es

\pf.  First of all, rather than working on the orbifold $\hat \M^d_\Dol$ and
restricting to $\hat P^d$, it is equivalent to work with
$\Ga$-equivariant objects on $\M^d_\Dol$ and restrict to $P^d$.

To show that $\hat B|_{P^d}$ is trivial, it suffices to show that the
projective bundle $\Psi|_{P^d}$ lifts to a $\Ga$-equivariant vector
bundle with trivial determinant.  We know from Lemma \re{cc} that
$\Psi|_{P^d}$ is the projectivization of 
a vector bundle,
but we need to show that $\Ga$ acts on this vector bundle.

Take a universal bundle over $\tilde J^d \times \tilde C$; 
since by Lemma \re{qq}
$$\tilde J^d \cong \frac{P^d \times J^0}{\Ga},$$ 
the pullback of this bundle to $P^d \times J^0 \times \tilde C$ has a
natural $\Ga$-action, which of course can be regarded as a $\tilde
\Ga$-action where the scalars $\C^\times$ act trivially.  By the
theorem of the cube \cite{mum} this pullback is equivariantly isomorphic to
$p_{13}^*\tilde \L \otimes p_{23}^* (1 \times \pi)^* \L^{-1}$, where
$p_{13}$ and $p_{23}$ are projections on the relevant factors, and
$\pi: \tilde C \to C$ is the spectral cover.  By Lemma \re{pp},
$\tilde \Ga$ acts on the second factor in this tensor product, with
the scalars acting with weight $-1$.  Hence it also acts on the first,
with the scalars acting with weight $1$.  Restricting to the basepoint
in $J^0$ gives us a $\tilde \Ga$-action on $\tilde \L \to P^d \times
\tilde C$, and hence on $\pi_* \tilde \L \to P^d \times C$.  Since
$\tilde \Ga_\c \cong \Ga \times \C^\times$ as mentioned before, this
produces a $\Ga$-action on $\pi_* \tilde \L |_{P^d \times \{ c \} }$.

Finally, as in the proof of Lemma \re{cc}, observe that $\det \pi_*
\tilde \L |_{P^d \times \{ c \} } \in \Pic^0_\Ga P^d = \Pic^0 \hat
P^d$.  So, by tensoring by an equivariant line bundle, the determinant
may be made equivariantly trivial.  \fp

Again we may examine the equivalence classes of
$\U{1}$-trivializations, which now form a torsor for $H^1(\hat
P^d,\U{1}) \cong P^0$.

\bs{Proposition} For any $d,e \in \Z$, there is a smooth isomorphism
of $P^0$-torsors
$$\Triv^{\U{1}}(\hat P^d,\hat B^e) \cong P^e.$$
\es

\pf.  First, it suffices as in the proof of Proposition \re{ee}
to take $e=1$.

Second, rather than working on $\hat \M^d_\Dol$, it is equivalent to
work $\Ga$-equivariantly on $\M^d_\Dol$.  For example, the torsor
$\Triv^{\Z_r}(\hat P^e,\hat B)$ parametrizing trivializations over
$\hat P^d$ can be identified with the torsor $\Triv_\Ga^{\Z_r}(P^e,B)$
parametrizing $\Ga$-equivariant trivializations over $P^d$.  As seen in
the proof of the previous lemma, such trivializations are provided
by $\tilde \Ga$-equivariant universal bundles $\tilde \L \to P^d
\times \tilde C$ where the scalars $\C^\times \subset \tilde \Ga$ act
with weight 1 and $\det \pi_*\tilde \L |_{P^d \times \{ c \} }$ is
trivial.  Indeed, the isomorphism classes of such universal bundles
form a torsor for $\Pic^0_\Ga P^d[r] = P^0[r]$, which must be
precisely $\Triv_\Ga^{\Z_r}(P^e,B)$.

Now follow the proof of Proposition \re{ee}: let $\mathcal T$ be the
torsor for $P^0 = \Pic^0_\Ga P^d$ parametrizing bundles $\tilde \L$
that satisfy all the conditions of the previous paragraph save that
$\det \pi_*\tilde \L |_{P^d \times \{ c \} }$ need only lie in
$\Pic^0_\Ga P^d$.  This contains the aforementioned $P^0[r]$-torsor
and hence must be isomorphic to $\Triv^{\U{1}}(\hat P^e,\hat B)$.

It remains only to identify ${\mathcal T}$ with
$P^1$.  First, recall that the $\tilde J^0$-torsor of all universal
bundles $\tilde \L \to \tilde J^d \times \tilde C$ with $\tilde
\L|_{\tilde J^d \times \{ \c \} } \in \Pic^0 \tilde J^d$ is isomorphic
to $\tilde J^1$, as shown in the proof of Proposition \re{ee}.

Then notice that there is an inclusion of ${\mathcal T}$ into this
torsor compatible with the inclusion $P^0 \subset \tilde J^0$.  It is
given simply by tensoring $\tilde \L \to P^d \times \tilde C$ by the
fixed $\tilde \Ga$-equivariant bundle $\pi^* \L^{-1} \to J^0 \times
\tilde C$, $\L$ being the universal bundle on $J^0 \times C$, to get a
$\Ga$-equivariant bundle over $P^d \times J^0 \times \tilde C$, which
descends to a universal bundle on the quotient $\tilde J^d \times
\tilde C$.

So ${\mathcal T}$ and $P^1$ are now both $P^0$-subtorsors of $\tilde
J^1$.  The quotient by either is the constant torsor $J^0$.  The image
of one in the quotient by the other therefore gives a morphism from
the base $U$ to $J^0$, the Jacobian of $C$.  But $U$ is a Zariski open
set in an affine space, so its only morphisms to an abelian variety
are constants.  Indeed, any nonconstant morphism would be nonconstant
on some line nontrivially intersecting $U$; the closure of the image
of this line would then be a rational curve in $J^0$, which doesn't
exist.

Hence ${\mathcal T}$ and $P^1$ are translates of one another in
$\tilde J^1$, so they are isomorphic.  \fp

We may summarize the results of this section as follows.

\bs{Theorem}
\label{ll}
For any $d, e \in \Z$, the moduli spaces $\M^d_\DR(\SL{r})$ and
$\M^e_\DR(\PGL{r})$, equipped with the flat unitary orbifold gerbes
$B^e$ and $\hat B^d$ respectively, are SYZ mirror partners.  \fp \es

\bit{Stringy mixed Hodge numbers}
\label{hodge}

Since the spaces we study are non-compact and singular, their ``Hodge
numbers'' must be interpreted in a generalized sense: as {\em stringy
  mixed Hodge numbers}.  Mixed Hodge numbers are alternating sums of
dimensions of the associated graded spaces in Deligne's mixed Hodge
structures on cohomology.  They are defined for any complex algebraic
variety, even incomplete or singular ones.  However, for the varieties
with orbifold (or more generally, Gorenstein) singularities arising in
string theory, mixed Hodge numbers are not the appropriate notion:
rather, we need a stringy version to take proper account of the
singularities.  For complete smooth varieties, these stringy mixed
Hodge numbers will coincide with the ordinary Hodge numbers.  More
generally, they coincide with the Hodge numbers of a crepant
resolution, if this exists \cite{craw}.

It is convenient to encode the mixed Hodge numbers as coefficients of
a polynomial: the so-called {\em $E$-polynomial}, or {\em virtual
Hodge polynomial}.   We will define a stringy $E$-polynomial in
terms of the ordinary one.

\noindent {\bf \boldmath The stringy $E$-polynomial.} 
The stringy $E$-polynomial is defined for any Gorenstein variety, but
it is expressed by a particularly simple formula in the case of a
quotient $M/\Ga$, where $M$ is a quasi-projective Calabi-Yau
$n$-manifold on which the finite group $\Ga$ acts preserving the
holomorphic $n$-form $\Om$.  We will treat this formula as a
definition, and present a generalization for $M$ equipped with a flat
unitary orbifold gerbe $B$.  The problems of how to interpret this
generalization in terms of smoothings of $M$, and how to extend it to
arbitrary Gorenstein varieties, are of the utmost interest, but we do
not pursue them here.

To any complex variety $X$, not necessarily smooth or projective,
Deligne \cite{del2,del3} has associated a canonical mixed Hodge structure
on the compactly supported cohomology $H^*_{\op{cpt}}(X,\C)$, and
hence, passing to the associated graded, complex vector spaces
$H^{p,q;k}(X)$.  These agree with $H^{p,q}(X)$ in the smooth
projective case, but in general they can be nonzero even when $p+q
\neq k$.  If a finite group $\Ga$ acts on $X$, it acts as well on each
$H^{p,q;k}(X)$; denote by $h^{p,q}(X)^\Ga$ the alternating sum over
$k$ of the dimensions of the $\Ga$-invariant subspaces.

Then define $E(X)^\Ga$ to be the polynomial in $u$ and $v$ given by
$$E(X)^\Ga = \sum_{p,q} h^{p,q}(X)^\Ga u^p v^q.$$
When $\Ga = 1$, this is the virtual Hodge polynomial $E(X)$ as defined
by, for example, Batyrev-Dais \cite{bd}.  A practical method of
determining $E(X)^\Ga$, which we adopt in the proof of Proposition
\re{ff}, is to regard $E(X)$ as a polynomial with coefficients in the
characters of $\Ga$, and then compute $E(X)^\Ga$ as the average value
over $\Ga$.

The beauty of $E(X)^\Ga$, like $E(X)$, is that it is additive for
disjoint unions and multiplicative for Zariski locally trivial
fibrations: the proofs given by Batyrev-Dais \cite[3.4, 3.7]{bd}, for
example, adapt without change to the equivariant case.  This allows us
to compute effectively in many cases even where we know nothing about
the mixed Hodge structures.

For $M$ as above, we may now define
the {\em stringy $E$-polynomial} to be
$$\Est(M/\Ga) = \sum_{[\ga]} E(M^\ga)^{C(\ga)}
(uv)^{F(\ga)}.$$ Here the sum runs over the conjugacy classes of
$\Ga$; $C(\ga)$ is the centralizer of $\ga$; $M^\ga$ is the subvariety
fixed by $\ga$; and $F(\ga)$ is an integer called the {\em fermionic
shift}, which is defined as follows.
The group element $\ga$ has finite order, so it acts on $TM|_{M^\ga}$
as a linear automorphism with eigenvalues $e^{2 \pi i w_1}, \dots,
e^{2 \pi i w_n}$, where each $w_j \in [0,1)$.  Let $F(\ga) = \sum
w_j$; this is an integer since, by hypothesis, $\ga$ acts trivially on
the canonical bundle.  (Purely for convenience of notation, we
have assumed that $F(\ga)$ is the same on all components of
$M^\ga$; otherwise we would have to write a further sum, over these
components, in the definition of $\Est$.)  

\noindent {\bf Turning on the \boldmath $B$-field.}  A twisted version
of this expression can be formulated in the following way.  Let $B$ be
an orbifold $\U{1}$-gerbe on $M/\Ga$, or equivalently, a
$\Ga$-equivariant $\U{1}$-gerbe on $M$.  Such a gadget induces a flat
$C(\ga)$-equivariant line bundle $L_{B,\ga}$ on the fixed-point set of
$\ga$.  Indeed, $B|_{M^\ga} \cong \ga^* B|_{M^\ga} = B|_{M^\ga}$,
where the isomorphism is given by the $\Ga$-action on $B$, and the
equality is because $\ga$ acts trivially on $M^\ga$.  This gives an
automorphism of $B$ restricted to $M^\ga$, and moreover, it is
$C(\ga)$-equivariant.  Any automorphism of a $\U{1}$-gerbe is given by
tensorization by a unique $\U{1}$-torsor, and this remains true
equivariantly.  Thus is determined a $C(\ga)$-equivariant
$\U{1}$-torsor on $M^\ga$, which is $L_{B,\ga}$.

We then propose the definition
\beq
\label{hh}
\Est^B (M / \Ga) = \sum_{[\ga]}
E(M^\ga;L_{B,\ga})^{C(\ga)}(uv)^{F(\ga)},
\eeq
where the $E$-polynomial is defined in terms of mixed Hodge numbers as
before, but on the cohomology with local coefficients in $L_{B,\ga}$.

Note that when $\ga=1$, the flat line bundle $L_{B,\ga}$ is
equivariantly trivial.
So we can regard the formula as saying
$$\Est^B (M / \Ga) = E(M)^\Ga + \cdots$$
where the dots denote the ``higher terms'' obtained from the fixed
points of $\ga \neq 1$.  In particular, viewing a smooth $M$ as
$M/\{1\}$, we find that $\Est^B(M) = E(M)$ for any flat gerbe $B$.
That is, the $B$-field affects the Hodge numbers only in the singular
case.

The case where $B$ is pulled back from a point is already nontrivial.
Indeed, we saw in \S\ref{syz} that $\Ga$-equivariant gerbes on a point
are classified up to isomorphism by the discrete torsion
$H^2(\Ga,\U{1})$.  For such a gerbe $B$, each bundle $L_{B,\ga}$ is
trivial and hence is determined by a $\U{1}$-representation of
$C(\ga)$.  This turns out to be $\delta \mapsto \nu(\delta,\ga) /
\nu(\ga,\delta)$, where $\nu$ is any group cocycle representing $B$.
The stringy $E$-polynomial therefore agrees in this case with the one
defined by Ruan \cite{ruan}.  But we will never use this fact.  Our
gerbes are not generally pulled back from a point, and in any case we
will construct the line bundles $L_{B,\ga}$ directly.  It does so
happen, though, that we get the same answer as we would from a certain
element of discrete torsion (cf.\ \cite{ht}).

\bit{The main conjecture}
\label{0}

Our purpose is to study the stringy mixed Hodge numbers of the moduli
spaces $M^d_\DR(\SL{r})$ and $M^d_\DR(\PGL{r})$.
We will assume, {\em now and henceforth}, that $r$ and $d$
are coprime.  Since $M^d_\DR(\SL{r})$ is $\Ga$-equivariantly
diffeomorphic to $M^d_\Dol(\SL{r})$, we may regard
$M^d_\DR(\SL{r})$ and $M^d_\DR(\PGL{r})$ as being equipped with the
flat unitary gerbes $B$ and $\hat B$ defined in \S\ref{triv}.  We then
conjecture the following.

\bs{Conjecture}
\label{jj}
For all $d,e \in \Z$ coprime to $r$,
$$\Est^{B^e}\Left(M^d_\DR(\SL{r})\Right) 
= \Est^{\hat B^d}\Left(M^e_\DR(\PGL{r})\Right).$$
\es

Since $M^d_\DR(\SL{r})$ is smooth, the left-hand side actually equals
$E\Left(M^d_\DR(\SL{r})\Right)$, which of course is independent of
$e$.

The rest of the paper is devoted to proving this conjecture for $r=2$
and $3$.  In fact much of what we prove is valid for general $r$.  
The broad outline of the argument is as follows.

First, we show in \S\ref{equal} that the stringy Hodge numbers of the
De Rham and Dolbeault spaces are the same.  Thereafter we may work
exclusively with the Dolbeault space, which has the advantage of
admitting a $\cx$-action.  So we wish to show
$$E\Left( \M^d_\Dol(\SL{r}) \Right) 
= \Est^{\hat B^d}\Left( \M^e_\Dol(\PGL{r})\Right).$$

In fact both sides are cumbersome to write down in full due to the
presence of a complicated ``leading term''
$E(\M^e_\Dol(\SL{r}))^\Ga$: the part invariant under the
$\Ga$-action.  But the remaining terms are more tractable.
So we will actually subtract it off and verify that 
\beq
\label{t}
E\Left(\M^d_\Dol(\SL{r})\Right) 
- E\Left(\M^e_\Dol(\SL{r})\Right)^\Ga =
\Est^{\hat B^d}\Left(\M^e_\Dol(\PGL{r})\Right) 
- E\Left(\M^e_\Dol(\SL{r})\Right)^\Ga.
\eeq

To compute the right-hand side, we need to know about the fixed-point
set of the $\Ga$-action.  This is described in \S\ref{I}, and the
computation is carried out for $r$ prime in \S\ref{II}.
To compute the left-hand side, we need to know about the fixed-point
set of the $\cx$-action.  This is described in \S\ref{IV}, and the
computation is carried through far enough to settle the
cases $r=2$ and $3$ in \S\ref{V}.

\bit{Equality of \boldmath $\Est$-polynomials of the Dolbeault and de
  Rham spaces} 
\label{equal}

For brevity, in this section $\M^d_\Hod(\SL{r})$ and
$\M^d_\Hod(\PGL{r})$ will be denoted simply by $\M_\Hod$ and $\hat
\M_\Hod$, respectively, and likewise for the Dolbeault and de Rham
spaces.

\bs{Lemma}
\label{u}
There exists a proper family $\overline{\M}_\Hod \to \C$
containing a divisor $X \times \C \to \C$ whose complement is
$\M_\Hod$.  It is a smoothly trivial family of orbifolds in the
sense that it is an orbifold, diffeomorphic to an orbifold times $\C$.
\es

\pf.  As seen in \S\ref{higgs}, $\cx$ acts on $\M_\Hod$ over the
action on the base $\C$ by scalar multiplication.  Let $\cx$ also act
on $\C^2$ by $t\cdot(x,y) = (tx,y)$.  Then $(x,y) \mapsto xy$ is a
$\cx$-equivariant map $\C^2 \to \C$.  Let $\M'$ be the base change of
$\M_\Hod$ given by pulling back by this map; then $\cx$ acts on $\M'$.
Regarded as a scheme over the second factor $\C$, the fiber of $\M'$
over $y \neq 0$ is $\M_\Hod$, but the fiber over $y=0$ is $\M_\Dol
\times \C$, with the diagonal action of $\cx$.

For any $p \in \M'$, the limit $\lim_{t \to 0} t \cdot p$ exists by
Corollary 10.5 of Simpson \cite{santa}.  Moreover, the fixed-point set
is $\M_\Dol^\cx \times \C \to \{ 0 \} \times \C$, which is proper over
$\C$ by Lemma 10.6 of Simpson \cite{santa}.  The hypotheses of Theorem
11.2 of Simpson \cite{santa} therefore hold, implying that the open
set $U \subset \M'$ of those $p \in \M'$ having no $\lim_{t \to
\infty} t \cdot p$ has a geometric quotient, proper and separated over
$\C$.  This open set is the complement of ${\mathcal N} \times \C \subset
\M_\Dol \times \C \to \{ 0 \} \times \C$, where $\mathcal N$ is the
so-called {\em nilpotent cone} in $\M_\Dol$, the zero fiber of the
Hitchin map.
The quotient $U / \cx$ is the disjoint union of two
pieces: an open set is the quotient of $\M' \sans (\M_\Dol \times \C)
\cong \M_\Hod \times \cx$, which of course is just $\M_\Hod$.  The
remainder is the quotient of $(\M_\Dol \sans {\mathcal N}) \times \C$,
which is of the form $X \times \C$, where $X$ is the geometric
quotient of $\M_\Dol \sans {\mathcal N}$.

Hence the quotient is a proper family of schemes over $\C$ whose
nonzero fiber is a compactification of $\M_\DR$ by adding $X$ as a
divisor at infinity, and whose zero fiber is a compactification of
$\M_\Dol$ by adding $X$ as a divisor at infinity.  In fact these
compactifications are precisely those constructed by Simpson
\cite{santa} and the first author \cite{hausel1}, respectively.  

Certainly $\overline \M_\Hod$ is an orbifold, as a geometric quotient of a
smooth variety by a $\cx$-action with finite stabilizers.
A neighborhood of any point in the zero fiber is diffeomorphic to a
trivial family of orbifolds: just note that $\cx$ acts trivially on
the base $\C$ and use the usual local model.  Then the standard
argument showing that a submersion of compact manifolds is locally
trivial applies in this orbifold situation: choose a Riemannian metric
and flow in a perpendicular direction to the projection.  So the
family is smoothly trivial in an analytical neighborhood of the zero
fiber.  But the $\cx$-action can be used to retract all of $\overline
\M_\Hod$ into this neighborhood. \fp

\bs{Theorem}
\label{v}
For $r$ and $d$ coprime,
$E(\M_\Dol) = E(\M_\DR)$.
\es

\pf.  The family constructed in the lemma above is a family of compact
``rational homology manifolds'' in the sense of Deligne
\cite[(8.2.4)]{del3}.  The mixed Hodge structures of the fibers are
therefore pure, that is, $H^{p,q;k} = 0$ unless $p+q=k$, and
Poincar\'e duality identifies the mixed Hodge structures on the
ordinary and compactly supported cohomology.  Because of the
topological triviality, the restriction from $\overline \M_\Hod$ to
any fiber is an isomorphism on cohomology, and hence an isomorphism of
mixed Hodge structures \cite[3.2.5]{del2}.  Hence the
mixed Hodge structures of $H^*_\cpt(\overline \M_\Dol)$ and
$H^*_\cpt(\overline \M_\DR)$ are isomorphic, and so $E(\overline
\M_\Dol) = E(\overline \M_\DR)$ .  But $\overline \M_\Dol$ is a
disjoint union $\M_\Dol \cup X$, while $\overline \M_\DR$ is a
disjoint union $\M_\DR \cup X$.  Since the $E$-polynomial is additive
under disjoint union, it follows that $E(\M_\Dol) = E(\M_\DR)$. \fp

\bs{Lemma}
\label{dd}
For any $\ga \in \Ga$, $\overline \M_\Hod^\ga$ is a smoothly trivial family of
orbifolds with $\overline \M_\Hod^\ga \cap (X \times \C) = X^\ga \times \C$.
\es

\pf.  The whole argument of Lemma \re{u} goes through provided that
$\M_\Hod^\ga$ is the geometric quotient by $\cx$ of $U^\ga$, where $U
\subset \M'$ is the open set in the proof of Lemma \re{u}.  In other words, we
want to know that $(U/\cx)^\ga = U^\ga/\cx$.  This means that if a
$\cx$-orbit is preserved by $\ga$, then it is fixed pointwise.  This
is obvious if the orbit does not lie over the $y$-axis in $\C^2$,
since $\Ga$ acts trivially on $\C^2$ while $\cx$ acts by $t \cdot
(x,y) = (tx,y)$.  On the other hand, the part of $U$ lying over any
point on the $y$-axis is $\M_\Dol \sans {\mathcal N}$, the complement
of the zero fiber of the Hitchin map.  But the Hitchin map $\mu: \M_\Dol
\to V_r$ takes the $\cx$-action on $\M_\Dol$ to a linear
action on the vector space $V_r$ with positive weights, while it
takes the $\ga$-action to the trivial action on $V_r$.  So the
only way for a $\cx$-orbit outside the zero fiber to be preserved by
$\ga$ is to be fixed pointwise.  \fp

\bs{Theorem} For any $e \in \Z$,
$$\Est^{\hat B^e}(\hat \M_\DR) = \Est^{\hat B^e} (\hat \M_\Dol).$$
\es

\pf.  Both sides are sums over $\ga \in \Ga$ by definition; it will be
shown that the terms agree, that is,
$$\Est(\hat M_\DR^\ga,L_{B^e,\ga}) = \Est(\hat
M_\Dol^\ga,L_{B^e,\ga}).$$
(The equality of the fermionic shifts is clear since
the representations of the finite group $\Ga$ are rigid.)

We wish to argue as in the proof of Theorem \re{v}, but first we need
to show that $L_{B^e,\ga} \to \M_{\Hod}^\ga$ extends over $\overline
\M_{\Hod}^\ga$ as a $\Ga$-equivariant flat line bundle.  Since
$\overline \M_{\Hod}$ is constructed in the proof of Lemma \re{u} as a
geometric quotient of the open set $U \subset M'$ described there, for
this it suffices to establish two statements: first, that the $\Ga$- and
$\cx$-actions on $L_{B^e,\ga}$ commute, and second, that the isotropy of
the $\cx$-action on $U$ acts trivially on $L_{B^e,\ga}$.

The first statement is easy: just notice that since $B^e$ is a
$\Z_r$-gerbe, $L_{B^e,\ga}$ has disconnected structure group $\Z_r$,
whereas $\cx$ is connected.  Since $1 \in \cx$ certainly commutes with
the $\Ga$-action, the whole of $\cx$ must.

As for the second statement, note that for any $p \in U$, the isotropy
group of the limit $\lim_{t \to 0} t \cdot p \in M'$ is $\cx$, and by
the same connectedness argument as in the previous paragraph this
isotropy group acts trivially on $L_{B^e,\ga}$.  Hence by continuity the
isotropy groups of $t\cdot p$, even though they may be disconnected,
must also act trivially.

Now proceed as in the proof of Theorem \re{v}, using the mixed Hodge
structure on cohomology with local coefficients \cite{ara,timm}.
Lemma \re{dd} guarantees that the same scheme $X^\ga$ gets added at
infinity to compactify both the De Rham and the Dolbeault fibers.  \fp

\bit{Fixed points of the \boldmath $\Sig$-action}
\label{I}

The action of $\Sig$ on the moduli space of stable bundles was
studied in a wonderful paper of Narasimhan and Ramanan \cite{nr}.  The
arguments in \S3 of their paper carry over without change to the space
of Higgs bundles.

Fix $\ga \in \Ga$ and let $m$ be its order.  Let $\pi: \tilde C \to C$
be the unbranched cyclic cover consisting of the $m$th roots of unity
in the total space of $L_\ga$.  A bundle on $C$ is equivalent to a
$\Z_m$-equivariant bundle on $\tilde C$, where $\Z_m$ is the Galois
group.

Let $(\bee,\bphi)$ be a universal Higgs bundle on $\tilde
M^d_\Dol(\GL{r/m}) \times \tilde C$, where the tilde denotes a moduli
space of bundles on $\tilde C$.  Then $\bphi$ induces a Higgs field
on $\pi_* \bee$; call it $\pi_* \bphi$, and regard $(\pi_* \bee, \pi_*
\bphi)$ as a family of Higgs bundles on $C$ parametrized by $\tilde
M^d_\Dol(\GL{r/m})$.  More precisely, note that as families of
$\Z_m$-equivariant bundles on $\tilde C$
$$\pi^*\pi_* \bee \cong \bigoplus_{i=1}^m (1 \times \xi^i)^* \bee,$$
where $\xi$, the standard generator of $\Z_m$, acts on the right-hand
side by cyclically permuting the factors; then the block-diagonal
Higgs field $\bigoplus_{i=1}^m (1 \times \xi^i)^* \bphi$ on the
right-hand side descends to the Higgs field we have called $\pi_*
\bphi$.  There is therefore an induced morphism $\tilde
M^d_\Dol(\GL{r/m}) \to M^d_\Dol(\GL{r})$.  Moreover, if $\delta
\in \Ga$ acts on $\tilde M^d_\Dol(\GL{r/m})$ by tensorization by $\pi^*
L_\delta$, then this morphism is $\Ga$-equivariant.

\bs{Proposition}
The action of $\Z_m$ on $\tilde \M^d_\Dol (\GL{r/m})$ is free, and the
morphism to $\M^d_\Dol(\GL{r})$ induced by $(\pi_* \bee, \pi_* \bphi)$
descends to a $\Ga$-equivariant regular embedding $$\tilde \M^d_\Dol
(\GL{r/m})/\Z_m \to \M^d_\Dol(\GL{r})$$ whose image is the fixed-point
set $M^d_\Dol(\GL{r})^\ga$.  
\es

\pf.  This proposition is analogous to Proposition 3.3 of
Narasimhan-Ramanan, and the proof is entirely similar.  The open set
$U$ that appears in their statements is unnecessary for us, since we
are assuming that $r$ and $d$ are coprime. \fp

Our next task is to ``fix the determinant,'' that is, pass from
structure group $\GL{r}$ to $\SL{r}$.  This requires some basic
facts about the Prym variety of an unbranched cyclic cover.  The
proofs of the following are pleasant exercises, and copious hints can
be found in Arbarello et al.  \cite[Appendix B2]{acgh}.

\begin{itemize}
  
\item Let $\tilde J^d = \Pic^d \tilde C$ and $J^d = \Pic^d C$.  Then
  the kernel of $\pi^*: J^0[m] \to \tilde J^0[m]$ is generated by
  $\ga$ [Exercise 14 in Arbarello et al.].
  
\item The kernel of the norm map $\Nm: \Pic \tilde C \to \Pic C$ has
  $m$ components [Exercise 19].  Call the identity component the {\em
  Prym variety} $P$.

\item The map $\Pic \tilde C \to \ker \Nm$ given by $L \mapsto L^{-1}
  \otimes \xi^* L$, where $\xi \in \Z_m$ is the standard generator,
  is surjective [Exercise 20].
  
\item For $\delta \in J^0[m]$, $\pi^* L_\delta$ is in the image of
  $\tilde J^d$ if and only if $\langle \ga, \delta \rangle = \xi^d$,
  where $\langle \,\, ,\, \rangle$ is the Weil pairing or intersection
  form on $J^0[m] = H_1(C,\Z_m)$ [Exercise 23].
  
\item For $L \in J^d$ with $(m,d) = 1$, the Galois group $\Z_m$ of
  $\tilde C \to C$ acts transitively on the set of components of
  $\Nm^{-1}(L)$ \cite[Proposition 3.5]{nr}.

\end{itemize}

There are natural splittings $\M^d_\Dol(\cx) = \Pic^d C \times
H^0(C,K)$ and $\tilde \M^d_\Dol(\cx) = \Pic^d \tilde C \times
H^0(\tilde C, K)$.  Define $\Pi: \tilde \M^d_\Dol(\cx) \to
\M^d_\Dol(\cx)$ to be $\det \pi_*$ on the first factor and the obvious
sum map on the second factor.  It is easy to see that $\det \pi_*$
equals $\Nm$ if $m$ is odd, and $\Nm$ composed with tensorization by
$L_\ga^{m/2}$ if $m$ is even: see Narasimhan-Ramanan for details.
Hence the fibers of $\Pi$ are torsors for $T^* \ker \Nm$ (over a
point).

\bs{Lemma}
The map induced by $(\pi_* \bee, \pi_* \bphi)$ lies in the following
commutative diagram:
$$\begin{array}{ccc}
\tilde \M^d_\Dol(\GL{r/m}) & \lrow & \M^d_\Dol(\GL{r}) \\
\down{\det} && \down{\det} \\
\tilde \M^d_\Dol(\cx) & \stackrel{\Pi}{\lrow} & \M^d_\Dol(\cx).
\end{array}$$
\es

\pf.  This is analogous to Lemma 3.4 of Narasimhan-Ramanan.  \fp

\bs{Corollary}
\label{ii}
If $m=r$, then the fixed-point set 
$\M^d_\Dol(\SL{r})^\ga$
is the quotient by $\Z_m$ of $\Pi^{-1}(L,0)$ for $L \in J^d$, which can
be identified with a connected component of $\Pi^{-1}(L,0)$, or with
the total space of the cotangent bundle of $\Nm^{-1}(L)$.  It is a
torsor for $T^*P$ (over a point).
\es

\pf. Since by definition $M^d_\Dol(\SL{r}) = \det^{-1}(L,0)$, this
follows immediately from the lemma and the facts preceding it. \fp

The identification in the corollary above is certainly convenient, but
it complicates the $\Ga$-action slightly.  Tensorization by a line
bundle of the form $\pi^* L_\delta$ may interchange the components of
the fiber of $\pi$, and we then have to act by an element of $\Z_m$ to
get back into our chosen one.  The following result clarifies how
this works.

\bs{Proposition}
\label{x}
{\rm (i)} Let $L \in J^d$ where $(m,d) = 1$, and let $q \in \Z$
satisfy $qd \equiv 1 \, (\mbox{\rm mod } m)$.  Then the action by
tensorization of $J^0[m]$ on $\Nm^{-1}(L)$ is transitive on the set of
components, and $\delta \in J^0[m]$ acts on the components in the same
way as $\langle \ga, \delta \rangle^q \in \Z_m$.  {\rm (ii)} The
Galois group $\Z_m$ acts on the Lie algebra of the Prym as $g-1$
copies of the regular representation of $\Z_m$ minus its trivial
factor.  \es

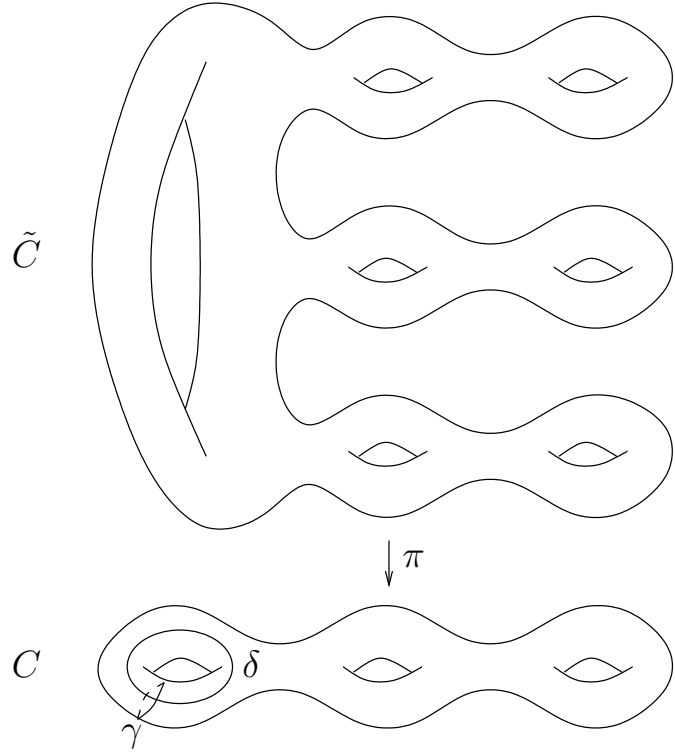
\pf.  The most enjoyable proof is topological.  Identify $J^0[m]$ with
$H_1(C,\Z_m)$; then $\Nm$ corresponds to the push-forward, and $\pi^*$
to the pullback of the Poincar\'e dual or inverse image, which we
denote by $\pi^{-1}$.  Since the intersection form is nondegenerate
and $\ga$ has order $m$, one can choose a set of generators for
$H_1(C,\Z_m)$ starting with $\ga$ so that the intersection form is
standard.  Consequently, there exists a handle presentation of $C$ so
that $\ga$ is represented by a loop around the first handle.  The
cover $\tilde C$ can then be depicted as in the diagram.

\begin{figure}
\begin{center}
\begin{picture}(0,0)%
\includegraphics{langlandsfigure.pstex}%
\end{picture}%
\setlength{\unitlength}{4144sp}%
\begingroup\makeatletter\ifx\SetFigFont\undefined
\def\x#1#2#3#4#5#6#7\relax{\def\x{#1#2#3#4#5#6}}%
\expandafter\x\fmtname xxxxxx\relax \def\y{splain}%
\ifx\x\y   
\gdef\SetFigFont#1#2#3{%
  \ifnum #1<17\tiny\else \ifnum #1<20\small\else
  \ifnum #1<24\normalsize\else \ifnum #1<29\large\else
  \ifnum #1<34\Large\else \ifnum #1<41\LARGE\else
     \huge\fi\fi\fi\fi\fi\fi
  \csname #3\endcsname}%
\else
\gdef\SetFigFont#1#2#3{\begingroup
  \count@#1\relax \ifnum 25<\count@\count@25\fi
  \def\x{\endgroup\@setsize\SetFigFont{#2pt}}%
  \expandafter\x
    \csname \romannumeral\the\count@ pt\expandafter\endcsname
    \csname @\romannumeral\the\count@ pt\endcsname
  \csname #3\endcsname}%
\fi
\fi\endgroup
\begin{picture}(3778,4499)(3421,-6910)
\put(4600,-6460){\makebox(0,0)[lb]{\smash{\SetFigFont{8}{9.6}{rm}{\color[rgb]{0,0,0}\large$\delta$}%
}}}
\put(3875,-6840){\makebox(0,0)[lb]{\smash{\SetFigFont{8}{9.6}{rm}{\color[rgb]{0,0,0}\large$\gamma$}%
}}}
\put(3221,-4003){\makebox(0,0)[lb]{\smash{\SetFigFont{8}{9.6}{rm}{\color[rgb]{0,0,0}\large$\tilde C$}%
}}}
\put(3221,-6460){\makebox(0,0)[lb]{\smash{\SetFigFont{8}{9.6}{rm}{\color[rgb]{0,0,0}\large$C$}%
}}}
\put(5556,-5797){\makebox(0,0)[lb]{\smash{\SetFigFont{8}{9.6}{rm}{\color[rgb]{0,0,0}\large$\pi$}%
}}}
\end{picture}

\end{center}
\end{figure}

The map whose image is the Prym is $L \mapsto L^{-1} \otimes \xi^* L$;
its restriction to $\tilde J^0[m] = H_1(\tilde C,\Z_m)$ is better
expressed in additive notation as $a \mapsto \xi^{-1} (a) - a$.  So we
need to find an element $\delta \in H_1(C,\Z_m)$ such that neither
$\pi^{-1} (\delta)$ nor any of its nonzero multiples are in the image of
this map.

The loop marked $\delta$ on the diagram clearly satisfies this
requirement.  On the other hand, the inverse images of all loops on
the last $g-1$ handles, equally clearly, are in the image of this map.
So the powers of $\delta$ act transitively on the set of components,
and the loops on the last $g-1$ handles act trivially.  Since $\langle
\ga, \delta \rangle = \xi \in \Z_m$ it now suffices for (i) to show
that $\delta$ acts as $\xi^q$.  Since the $\Ga$- and $\Z_m$-actions
clearly commute, $\delta$ must act as some power of $\xi$.  But
according to the penultimate basic fact, there exists $M \in \tilde
J^1$ such that $M^{-1} \otimes \xi^{-1}M \cong \pi^* L_\delta$, and
hence $\xi^* M^d \cong \pi^* L_{d\delta} \otimes M^d$.  Tensor $M$
by a line bundle pulled back from $J^0$ so as to arrange that $M^d \in
\Nm^{-1}(L)$.  Then the actions of $\xi$ and $d\delta$ agree on the
component containing $M^d$, and hence on all components.  Therefore
the same is true of $\xi^q$ and $\delta$.  This proves (i).

The Lie algebra of the Prym can be identified with $H^1(\tilde C,\R)
/ H^1(C,\R)$.  This is spanned by the loops on the last $m(g-1)$
handles of $\tilde C$, modulo the inverse images of the loops from the
last $g-1$ handles of $C$, and (ii) is now clear.  \fp

\bs{Corollary}
\label{gg}
Let $L$ and $q$ be as above, and identify $M^d_\Dol(\SL{r})^\ga$ with
the cotangent bundle of a connected component of $\Nm^{-1}(L)$.
Then the $\Ga$-action on $M^d_\Dol(\SL{r})^\ga$ is induced by the
following action on that connected component: $\delta \in \Ga$ acts by
tensorization by $\pi^* L_\delta$ followed by the action of
$\langle \ga,\delta \rangle^{-q} \in \Z_m$.  
This acts on $H^1(M^d_\Dol(\SL{r})^\ga, \R)$ as stated in {\rm (ii)}
above. \fp
\es

\bit{Calculation for \boldmath $\PGL{r}$}
\label{II}

Suppose that $r$ is prime and that $d$ and $e$ are coprime to $r$.
Then we can work out the right-hand side of \re{t} completely.  By
abuse of notation we refer henceforth to the fixed-point set
$M^d_\Dol(\SL{r})^\ga$ as $T^*P$, although it is really a torsor for
$T^*P$ over a point.

\bs{Proposition}
The $\Ga$-equivariant flat line bundle $L_{B,\ga} \to T^* P$ is trivial,
and the $\Ga$-action is given by the character $\delta \mapsto \langle
\ga,\delta \rangle^{-q}$, where $\Z_r$ is identified as usual with the
complex $r$th roots of unity.  
\es

\pf.  Since we are studying a flat line bundle, instead of working
with $T^* P$ we may work just with the zero section.  This is
convenient, since the Higgs field vanishes there, so we may forget it
and think of the universal family as merely a bundle.

We abusively call the zero section $P$, but it is really a component
of $\Nm^{-1}(L) \subset \tilde J^d$, which is a torsor for $P$.
According to Corollary \re{gg}, $\delta \in \Ga$ acts on it by
tensorization by $\pi^* L_\delta$ followed by the action of $\langle
\ga,\delta \rangle^{-q} \in \Z_m$.  We now explain how to lift both of
these actions to (projective) actions on the universal bundle.

First, let $\xi = e^{2\pi i/m} \in \Z_m$ act as an element of the
Galois group, both on $\Nm^{-1}(L) \subset \tilde J^d$ and on $\tilde
C$ itself.  Take a universal line bundle $\L \to \tilde J^d \times
\tilde C$; then by the universal property, $(\xi \times \xi)^* \L
\cong Q \otimes \L$ for some line bundle $Q \to \tilde J^d$.

The push-forward $\pi_* \L$
is the desired universal bundle.
Since $\pi: \tilde C \to C$ is a Galois cover,
there is an isomorphism of $\Z_m$-equivariant bundles
$$\pi^* \pi_* \L \cong \bigoplus_{i=1}^m (1 \times \xi^i)^* \L, $$
where on the right-hand side the factors are cyclically permuted by
the action.

Hence there are isomorphisms of $\Z_m$-equivariant bundles
\beqas
(\xi\times 1)^* \pi^*\pi_* \L 
& \cong &
\bigoplus_{i=1}^m (\xi \times \xi^i)^* \L \\
& \cong &
\bigoplus_{i=1}^m (1 \times \xi^{i-1})^* (\xi \times \xi)^* \L \\
& \cong &
\bigoplus_{j=1}^m (1 \times \xi^j)^* (Q \otimes \L),  \\
& \cong &
Q \otimes \pi^*\pi_* \L, 
\eeqas
where the penultimate step makes the change of variables $j=i-1$.
This descends to the desired isomorphism $(\xi\times 1)^* \pi_* \L \cong
Q \otimes \pi_* \L$.  

Restricting to $P \times \{c \}$, we find that the projective bundle
$\Psi|_P$ involved in the definition of the gerbe $B|_P$ is in fact the
projectivization of the vector bundle 
$$V = \bigoplus_{y \in \pi^{-1}(c)} \L_y$$ 
where $\L_y = \L|_{P \times \{ y \}}$, and that the projective
action of $\xi$ cyclically permutes the summands.

Second, let $T_\delta: \tilde J^d \to \tilde J^d$ denote tensorization by
$\pi^* L_\delta$, which preserves $\Nm^{-1}(L)$.  If $\L$ is chosen to
be trivial over a basepoint in $\tilde C$, then by the universal
property, $(T_\delta \times 1)^* \L \cong \pi^* L_\delta
\otimes \L$.  Hence there are isomorphisms of $\Ga$-equivariant bundles
\beqas 
(T_\delta \times 1)^* \pi^*\pi_* \L 
& \cong &
\bigoplus_{i=1}^m (T_\delta \times \xi^i)^* \L \\
& \cong &
\bigoplus_{i=1}^m (1 \times \xi^i)^* (\pi^* L_\delta \otimes \L)\\
& \cong &
\pi^* L_\delta \otimes \bigoplus_{i=1}^m (1 \times \xi^i)^* \L, \\
& \cong & 
\pi^* L_\delta \otimes\pi^*\pi_* \L
\eeqas
descending to an isomorphism $T_\delta^* \pi_* \L \cong
L_\delta \otimes \pi_* \L$.  

Restricting again to $P \times \{c \}$, we find that $T_\delta$ acts
on $V$ as an isomorphism on each summand.  For example, if
$\delta=\ga$, then $\pi^* L_\ga$ is the trivial bundle, with $\xi$
acting by multiplication by $e^{2\pi i /m}$; so for each $y \in
\pi^{-1}(c)$, the isomorphism $\L_{\xi\cdot y} \to \L_{\xi \cdot y}$
is $e^{-2\pi i / m}$ times the isomorphism $\L_y \to \L_y$.  In
particular, the automorphism of $\Psi|_P$ induced by the action of
$\ga$ lifts to an automorphism of $V$: in other words, it takes this
lifting to an isomorphic lifting.  This means that the flat line bundle
$L_{B,\ga}$ defined by the automorphism of the gerbe of liftings
$B|_P$ is trivial.

However, the action of $\Ga$ on this flat line bundle is not trivial.
Indeed, the action of $\delta$ as described in Corollary \re{gg} lifts to
a projective action on $\pi_* \L$ via the isomorphisms above.  This
leads to a diagram
$$ 
\begin{array}{ccc}
V & \stackrel{\ga}{\lrow} 
& L_{B,\ga} \otimes V \\
\leftdown{\delta} && \down{\delta} \\
V & \stackrel{\ga}{\lrow} 
& L_{B,\ga} \otimes V
\end{array}
$$
but, to make the diagram commutative, we must multiply $L_{B,\ga}$
by a scalar factor.  Since $\delta$ cyclically permutes the summands
as $\langle \ga, \delta \rangle^q \in \Z_m$ and $\ga$ acts on each
successive summand as $e^{-2\pi i/m}$ times the previous one, this
factor is $\langle \ga, \delta \rangle^{-q}$, as desired.  \fp

\bs{Proposition}
\label{ff}
When $r$ is prime, the 
right-hand side of \re{t} equals
$$\ratio{1}{r}(r^{2g}-1)(uv)^{(r^2-1)(g-1)}
\left(
(1-u)^{(r-1)(g-1)}(1-v)^{(r-1)(g-1)}
- \left( \frac{(1-u^r)(1-v^r)}{(1-u)(1-v)} \right)^{g-1}
\vphantom{\Left( \Right)^{g-1}}
\right).$$
\es

\pf.  The definition \re{hh} of stringy Hodge numbers calls for adding up a
contribution from the fixed-point set of each nontrivial $\ga \in \Ga$.
As seen in Corollary \re{ii}, this fixed-point set is (a torsor for)
$T^*P$.  The compactly supported cohomology of $T^*P \cong
\C^{(r-1)(g-1)} \times P$ splits according to the K\"unneth formula,
and that of the first factor is of course $\Ga$-invariant, so \beqas
E(T^*P,L_{B,\ga})^\Ga & = &
E(\C^{(r-1)(g-1)}) E(P,L_{B,\ga})^\Ga \\
& = & (uv)^{(r-1)(g-1)} E(P,L_{B,\ga})^\Ga.  \eeqas

To evaluate the right-hand side, note that by Corollary \re{gg}, any
$\delta \in \Ga$ acts on $H^1(P,\R)$, and hence on $H^{0,1}(P)$, with
eigenvalues $\langle \ga, \delta \rangle^k$, for $k = 1$ to $r-1$,
each repeated $g-1$ times.  Since
$$H^k(P,\C) = \La^k\Left( H^{0,1}(P) \oplus
H^{1,0}(P) \Right),$$ 
as a polynomial with coefficients in the characters of $\Ga$,
$$E(P) = 
\left(
\prod_{i=1}^{r-1}(1-\rho^i u)(1-\rho^i v)
\right)^{g-1}$$
where $\rho(\delta) = \langle \ga, \delta \rangle$, 
and $E(P,L_{B,\ga}) = \rho^{-e} E(P).$  
The invariant part is
the average value:
\beqas
E(P,L_{B,\ga})^\Ga 
&=&
\frac{1}{|\Ga|} \sum_{\delta \in \Ga} 
\rho^{-1}(\delta)
\left(
\prod_{i=1}^{r-1}(1-\rho^i(\delta) u)(1-\rho^i(\delta) v)
\right)^{g-1} \\
&=&
\ratio{1}{r} \sum_{i=0}^{r-1}
\xi^{-ei}
\left(
\prod_{i=1}^{r-1}
(1-\xi^i u)(1-\xi^i v)
\right)^{g-1} \\
&=&
\ratio{1}{r}
\left(
(1-u)^{(r-1)(g-1)}(1-v)^{(r-1)(g-1)}
- \left( \frac{(1-u^r)(1-v^r)}{(1-u)(1-v)} \right)^{g-1}
\vphantom{\Left( \Right)^{g-1}}
\right),
\vphantom{\frac{1}{1}}
\eeqas
where $\xi = e^{2\pi i/r}$.  

To compute the fermionic shift, note that $\ga$ acts with nontrivial
weights on the normal bundle to $T^*P$ in $M^d_\Dol(\SL{r})$.  The
action of $\ga$ preserves the holomorphic symplectic structure, since
on the dense open set in $\M^d_\Dol(\SL{r})$ isomorphic to the
cotangent bundle to the moduli space of stable bundles it corresponds
to the tautological symplectic structure.  Hence every eigenvalue
$e^{2\pi i \alpha}$ is accompanied by an eigenvalue $e^{2\pi
i(1-\alpha)}$, so the fermionic shift is half the rank of the normal
bundle, namely $r(r-1)(g-1)$.  Summing over the $r^{2g}-1$ identical
terms yields the grand total in the statement.  \fp

\bit{Fixed points of the \boldmath $\cx$-action}
\label{IV}

The Betti numbers of $\M^d_\Dol(\SL{r})$ are computed by Hitchin
\cite{lms} and Gothen \cite{goth} for ranks 2 and 3 respectively, and
the $E$-polynomials can be calculated in the same way.  But the complete
formula is complicated and unilluminating.  All we want to
know, as was explained in \S\ref{0}, is the Hodge polynomial of what
we like to call the {\em variant} cohomology: the part not invariant
under the action of $\Ga$.  This is given by the left-hand side of
\re{t}.

For convenience, in this section denote $\M^d_\Dol(\SL{r})$ simply by
$M_\Dol$.  To describe its variant cohomology, we shall consider the
action of the multiplicative group $\cx$ on $M_\Dol$ given by $\la
\cdot (E, \phi) = (E, \la \phi)$.  This commutes with the
$\Ga$-action.  Let $\mc F$ be the fixed-point set.

\bs{Proposition}
\label{e}
As polynomials with coefficients in the characters of $\Ga$,
$$E(M_\Dol) = (uv)^{\dim M_\Dol/2}E({\mc F}).$$
\es

\pf. The $\gli$-action satisfies the property that for all $x \in M$,
$\la \cdot x$ has a limit as $\la \to 0$.  This follows directly from
the properness of the Hitchin map, since it takes this
$\gli$-action to a linear action on $V$ with positive weights.  

Now there is an algebraic version of the Morse stratification called the
Bia\l ynicki-Birula stratification \cite{bb,flip}.  Indeed, in this
case it is nothing but the Morse stratification for the moment
map for the action of $\U{1} \subset \gli$. 

It implies that $M_\Dol$ is a $\Sig$-invariant union of Zariski
locally trivial fiber bundles whose fibers are affine spaces, whose
bases are the components of $\mc F$, and whose projections are given
by $x \mapsto \lim_{\la \to 0}\la \cdot x$.  Not only that, the
dimension of the affine space is always $\dim M_\Dol/2$.  One could
prove this directly by looking at the $\gli$-action on the deformation
space \cite{higgs1}.  A lazier proof, however, is just to quote
Ginzburg's result \cite{ginz} that the downward flow from each
critical set is Lagrangian, and hence has dimension equal to
half that of $M_\Dol$.  The same is therefore true of the upward flow
from each point in the critical set, which is the affine space.

The desired formula follows from the additivity of the $E$-polynomial
for disjoint unions and its multiplicativity for Zariski
locally trivial fibrations \cite[3.4, 3.7]{bd}. \fp

\bs{Lemma (Simpson)}
\label{z}
If $(E,\phi) \in {\mc F}$, then there exists a decomposition $E =
\bigoplus E_i$ with $\phi(E_i) \subset K \otimes E_{i+1}$.  Moreover,
the ranks and degrees of the $E_i$ are locally constant on $\mc F$.
\es

\pf.  Fix $t \in \gli$ which is not a root of unity.  If $(E,\phi)$ is
to be in $\mc F$, then there must be an isomorphism $f: E \to E$ such
that $f \phi = t \phi f$.  Such an $f$ is unique up to a scalar, since
two such maps $f$ and $f'$ give rise to an automorphism $f^{-1} f'$ of
the stable pair $(E,\phi)$, which must be a scalar.  The roots of
the characteristic polynomial form an $r$-fold cover of $C$ in $C
\times \C$, so they and their multiplicities are constant on $C$.
This gives a decomposition of $E$ into generalized eigenspaces
$E_\la$, the kernels of $(f-\la)^r$.  These $E_\la$ constitute
subbundles of the universal bundle $\bee$ restricted to each connected
component of $\mc F$.  For locally on $\mc F$, $f$ extends to an
automorphism of $\bee|_{{\mc F} \times C}$; indeed, the
hypercohomology $\Hyp^0$ of the two-term complex $\End E \to K \otimes
\End E$ on $C$ with $f \mapsto f \phi - t \phi f$ is one-dimensional
and generated by the $f$ mentioned above, so the hyper-direct image on
$\mc F$ $(\R^0 \pi)_* (\End \bee \to K \otimes \End \bee)$ is locally
free of rank 1.  Hence the ranks and degrees of the $E_\la$ are
locally constant on $\mc F$.

Now $(f-t \la)^r \phi = t^r \phi (f-\la)^r$, so $\phi$ maps the
$\la$-generalized eigenspace $E_\la$ to the $ \la$-generalized
eigenspace $E_{t \la}$.  Since $t$ is not a root of unity, the
eigenvalues break up into finite strings $\la, t\la, \dots, t^k \la$,
but as stable Higgs bundles are irreducible, there is only one such
string.  \fp

It will be convenient to refer to the finite sequence $(\rk E_1, \rk
E_2, \dots)$ as the {\em type} of the component of $\mc F$ containing
$(E,\phi)$.  One possibility is the type $(r)$ consisting of a single
number only.  This means that the Higgs field vanishes, so the
corresponding component is simply the moduli space of stable vector
bundles of rank $r$ and fixed determinant.

\bit{Calculation for \boldmath $\SL{r}$}
\label{V}

Now suppose once again that $r$ is prime.  We will calculate the
contribution to the variant cohomology of the fixed points of type
$(1,1,\dots,1)$.  Then we will show that in ranks 2 and 3, these are
the only nonzero contributors.

\bs{Proposition}
\label{kk}
When $r$ is prime, the fixed components of type $(1,1,\dots,1)$
contribute to the variant part of $E\Left(\M^d_\Dol(\SL{r})\Right)$
exactly the polynomial given in Proposition \re{ff}.  \es

\pf.  A Higgs bundle of type $(1,1,\dots,1)$ has the form $E =
\bigoplus_{i=1}^r L_i $ with $\phi_i:L_i \to L_{i+1} \otimes K$.  We
assume that the determinant is fixed to be $\O(dc)$ where $c$ is our
basepoint, so $\prod L_i \cong \O(dc)$.  Let $D_i$ be the divisor of
zeroes of $\phi_i$ and $M_i = \O(D_i)$.  Then $M_i \cong
L_i^{-1}L^{\phantom{-1}}_{i+1}K$ and so 
\beq 
\label{aa}
\prod_{i=1}^{r-1}M_i^i \cong L_r^r K^{r(r-1)/2}(dc).  
\eeq 
Denote
$l_i = \deg L_i$ and $m_i = \deg M_i$; then $m_i = l_i - l_{i+1} + 2g
-2$, and 
\beq 
\label{bb} 
\sum_{i=1}^{r-1} i m_i \equiv d \,\,\,\,
\mbox{(mod $r$)}.  
\eeq
By the way, this last constraint is accidentally overlooked in the
paper of Gothen \cite{goth}, leading to some incorrect formulas.

A $\cx$-invariant Higgs bundle is unstable if and only if it is
destabilized by a $\cx$-invariant Higgs subbundle: indeed, this
follows immediately from the uniqueness of the Harder-Narasimhan
stratification for Higgs bundles.  Since the only such
subbundles are of the form $\bigoplus_{i=k}^r L_i$, stability is
equivalent to
$$\frac{l_k+l_{k+1} + \cdots + l_r}{r-k+1} < \frac{d}{r}$$
for each $k$.
It is a simple exercise to show that these inequalities are satisfied
if $0 \leq m_i \leq 2g-2$.  These are the only values that will
contribute to the variant cohomology.

Given effective divisors $D_i$ whose degrees satisfy \re{bb}, all that
is needed to construct a Higgs bundle of the type described above is
a choice of $L_r$, which by \re{aa} is determined up to multiplication
by an $r$th root of unity, that is, an element of $\Ga$.
Consequently, each type $(1,1,\dots,1)$ component of the fixed-point
set is a fibered product of the form
$$N_{m_1,\dots,m_r} = \left( \prod_{i=1}^{r-1} S^{m_i} C \right)
\times_{\Pic^{\sum i m_i} C} \Pic^{l_r} C,$$
where the morphism from $\prod_{i=1}^{r-1} S^{m_i} C$
is $(D_i) \mapsto \O(\sum i D_i)$ and the morphism from $\Pic^{l_r} C$ is
$L \mapsto L^r K^{r(r-1)/2} (dc)$.

The terms in the Hodge decomposition, as representations of $\Ga$, can
be computed by pushing forward to $\prod_{i=1}^{r-1} S^{m_i} C$ first,
as in Hitchin \cite{lms} and Gothen \cite{goth}.  It turns out that 
\beqas
H^*(N_{m_1,\dots,m_r}, \C) 
&=& 
\bigoplus_{\ga \in \Ga}
H^*\Left(\prod_i S^{m_i} C,\bigotimes_i \pi_i^* \L^i_\ga\Right) 
\\ 
&=&
\bigoplus_{\ga \in \Ga} \bigotimes_i \La^{m_i}H^1(C,L_\ga^i), 
\eeqas
where the right-hand side denotes cohomology with local coefficients,
and $\L_\ga \to S^{m_1} C$ is the flat line bundle obtained either by
symmetrizing $L_\ga \to C$ or by pulling back the corresponding flat
line bundle over $\Pic^0 C$ via the Abel-Jacobi map.  The variant part
consists of the terms where $\ga \neq 1$.  Since $r$ is prime,
$L_\ga^i$ is then a nontrivial flat bundle for each $i$ from $1$ to
$r-1$, and hence $H^1(C,L_\ga^i)$ has Hodge type $(g-1,g-1)$; indeed,
its $(0,1)$ part can be identified with the Dolbeault cohomology of
$L_\ga^i$ on $C$.

The contribution of $N_{m_1,\dots,m_r}$ to the variant part of the
$E$-polynomial is therefore
$$(r^{2g}-1)(uv)^{(r^2-1)(g-1)}
\Coeff_{\prod t_i^{m_i}} \left( \prod_i (1-t_iu)(1-t_iv)
\right)^{g-1}.$$
Here the factor of $r^{2g}-1$ is the number of nontrivial group
elements in $\Ga$, and the power of $uv$ is the contribution of the
normal bundle, as described in \re{e}. 
To sum $m_i$ from $0$ to $2g-2$ subject to the constraint \re{bb}, let
$\xi = e^{2\pi i/r}$ and take the average value of this multiplied by
$\xi^{-d}$: that is,
\beqas
\lefteqn{\ratio{1}{r} (r^{2g}-1) (uv)^{(r^2-1)(g-1)}
\sum_{j=1}^r \xi^{-jd} 
\left( \prod_{i=1}^{r-1}(1-\xi^{ij}u)(1-\xi^{ij}v) \right)^{g-1}} \\
&=& \ratio{1}{r} (r^{2g}-1) (uv)^{(r^2-1)(g-1)}
\left( 
\left((1-u)(1-v) \right)^{(r-1)(g-1)}
- \left( \frac{(1-u^r)(1-v^r)}{(1-u)(1-v)} \right)^{g-1}
\right). 
\eeqas
This is indeed the polynomial given in Proposition \re{ff}.
\fp

\bs{Lemma}
\label{f}
In any rank $r$, the fixed component of type $(r)$ has no variant cohomology.
\es

\pf. This fixed component is the moduli space of stable
vector bundles of rank r and determinant $\O(dc)$, $c \in C$ being our
chosen basepoint.  So the desired fact is exactly Theorem 1 of
Harder-Narasimhan \cite{hn}, cf.\ also Newstead \cite{new} and
Atiyah-Bott \cite{ab}.  \fp

\bs{Lemma}
In rank $3$, the fixed components of type $(1,2)$ and $(2,1)$ have no
variant cohomology. 
\es

\pf.  Gothen \cite{goth} shows that each such fixed component is a
smooth family over $\Pic^0 C$ whose fiber is the moduli space of
stable rank 2 Bradlow pairs with a certain fixed determinant and a
fixed Bradlow parameter $\tau$.  As such, this family can be obtained,
starting from a projective bundle over $\Pic^0 C$, by a sequence of
blow-ups and blow-downs whose centers are projective bundles over
$\Pic^0 C$ times symmetric products of $C$, in the manner prescribed
by the second author \cite{pair}.  Gothen explains how all the spaces in
this sequence can be regarded as parametrizing families of (not
necessarily stable) rank 3 Higgs bundles, and it follows that $\Ga$
acts on all the spaces, compatibly with all the morphisms between
them.  Furthermore, it acts on the projective bundles by bundle maps,
and on their bases by translation of the factor $\Pic^0 C$: this is
readily apparent from Gothen's description.  Consequently, it acts
trivially on the cohomology of each projective bundle in the sequence.
The standard description of the cohomology of a blow-up (see e.g.\
Griffiths-Harris \cite[p.\ 605]{gh}) implies that, when a finite group
$\Ga$ acts on a smooth $X$ preserving a smooth $Y \subset X$, it acts
trivially on the cohomology of the blow-up along $Y$ if and only if it
acts trivially on the cohomology of $X$ and of $Y$.  Hence $\Ga$ acts
trivially on the cohomology of the fixed component.  \fp

\bs{Theorem}
\label{mm}
Conjecture \re{jj} holds true for $r=2$ and $3$.
\es

\pf. First of all, when $r = 2$ or $3$, and $d$ and $e$ are both
coprime to $r$, there is of course an isomorphism $\M^d_\Dol(\SL{r})
\cong \M^e_\Dol(\SL{r})$: it is given simply by dualizing (in the case
when $r=3$ and $d \not \equiv e$ mod 3) and tensorizing by a line
bundle of the appropriate degree.  Hence we may substitute $d$ for $e$
on the left-hand side of \re{t} with impunity.  It is then just a
question of studying the variant part of the $E$-polynomial for
$\M^d_\Dol(\SL{r})$.

The contribution of the $(1,1)$ or $(1,1,1)$ components to this
variant part, given by Proposition \re{kk}, agrees with the
calculation for $\PGL{r}$ given in Proposition \re{ff}.  The
contributions of the remaining components vanish by the two lemmas
above. \fp

More generally, for any prime $r$, Conjecture \re{jj} would follow from
two further conjectures.  First, that $E\Left(\M^d_\Dol(\SL{r})\Right)
= E\Left(\M^e_\Dol(\SL{r})\Right)$ for all $d$ and $e$ coprime to $r$.
Second, that {\em no\/} fixed component besides those of type
$(1,1,\dots,1)$ contributes to the variant cohomology.  We hope to
return to these conjectures in the future.

\end{document}